\magnification = 1200
\baselineskip = 15 pt

\input psfig.sty

\def \sect#1{\bigskip \noindent {\bf #1} \medskip}
\def \subsect#1{\medskip \noindent{\it #1} \medskip}
\def \th#1#2{\medskip \noindent {\bf Theorem #1.}   \it #2 \rm}
\def \prop#1#2{\medskip \noindent {\bf Proposition #1.}   \it #2 \rm}
\def \cor#1#2{\medskip \noindent {\bf Corollary #1.}   \it #2 \rm}
\def \pf {\noindent  {\it Proof}.\quad }
\def \lem#1#2{\medskip \noindent {\bf Lemma #1.}   \it #2 \rm}

\centerline{\bf Minimizing the Probability of Lifetime Ruin under Borrowing Constraints} \bigskip \bigskip

\centerline{Version: 12 October 2006}  \bigskip

\noindent Erhan Bayraktar  \hfil \break
\noindent Department of Mathematics, University of Michigan  \hfil \break

\noindent Virginia R. Young  \hfil \break
\noindent Department of Mathematics, University of Michigan \hfil \break
%\noindent 530 Church Street \hfil \break
%\noindent Ann Arbor, Michigan, 48109 \hfil \break
%\noindent Phone: 734-764-7227 \hfil \break
%\noindent Email: vryoung@umich.edu \hfil \break

\noindent{\bf Abstract:} We determine the optimal investment strategy of an individual who targets a given rate of consumption and who seeks to minimize the probability of going bankrupt before she dies, also known as {\it lifetime ruin}.  We impose two types of borrowing constraints:  First, we do not allow the individual to borrow money to invest in the risky asset nor to sell the risky asset short.  However, the latter is not a real restriction because in the unconstrained case, the individual does not sell the risky asset short.  Second, we allow the individual to borrow money but only at a rate that is higher than the rate earned on the riskless asset.

We consider two forms of the consumption function: (1) The individual consumes at a constant (real) dollar rate, and (2) the individual consumes a constant proportion of her wealth. The first is arguably more realistic, but the second is closely connected with Merton's model of optimal consumption and investment under power utility.  We demonstrate that connection in this paper, as well as include a numerical example to illustrate our results.

\medskip


\medskip


\vfill
\eject

\centerline{\bf Minimizing the Probability of Lifetime Ruin under Borrowing Constraints}

\sect{1. Introduction}

There is growing concern among Americans about financial ruin during retirement (Parikh, 2003).  These concerns are justified because a significant financial crisis is looming; it is projected that retired Americans' living expenses will exceed their financial resources by \$400 billion over the ten-year period 2020-2030 (VenDerhei and Copeland, 2003).  This shortfall is driven by demographic trends, the increased longevity of our aging population, changes in Social Security, inadequate private retirement savings, and the continuing trend toward defined contribution plans, such as 401(k)s, under which the individual - not the employer - assumes all investment and longevity risk.

We employ techniques of stochastic optimal control to study the problem of how an individual should invest her wealth in a risky financial market in order to minimize the probability that she outlives her wealth, also known as the probability of {\it lifetime ruin} (Milevsky and Robinson, 2000). Specifically, we determine the optimal investment strategy of an individual who targets a given rate of consumption and who seeks to minimize the probability of lifetime ruin. As employers shift from defined benefit plans to defined contribution plans, the problem of outliving one's wealth becomes relevant to more retirees and to the actuaries and financial planners who advise them.  A recent issue of {\it The Actuary} (Parikh, 2003) points out that ``according to the Employee Benefit Research Institute, in 1974, fifty-six percent of retirement income was coming from guaranteed sources. \dots Cerulli Associates projects that will drop to twenty-four percent by 2030.''

When an individual seeks to find an optimal investment policy, the resulting optimal policy depends on her optimization criterion.  The most common optimization criterion encountered in the finance literature is to  maximize one's expected discounted utility of consumption and bequest.  In the 1970s, Merton (1992) began study of this problem, and many others  continued his work; see, for example, Karatzas and Shreve (1998, Chapter 3) and the discussion at the end of that chapter for further references.  Notable extensions include the work of Davis and Norman (1990), Zariphopoulou (1992), and Shreve and Soner (1994) on portfolio selection with transaction costs; Duffie and Zariphopoulou (1993), Duffie et al. (1997), and Koo (1998) on optimal investment and consumption strategies to maximize expected utility of consumption and bequest in the presence of stochastic income; and Fleming and Zariphopoulou (1991) and Zariphopoulou (1994) on optimal investment under borrowing and trading constraints, respectively.  See Zariphopoulou (1999, 2001) for helpful summaries of the work to date in this area.

The goal of maximizing expected discounted utility of consumption and bequest may be difficult to implement because it depends on a subjective utility function for consumption and bequest.  Minimizing the probability of running out of money before dying might prove easier for individuals to understand because that criterion is arguably more objective.  Indeed, there is something intuitively appealing about minimizing the probability of shortfall that lends itself to asset allocation advice. In fact, the Nobel laureate William Sharpe founded a financial services advisory firm that is largely based on using probabilities to provide investment advice.

Other researchers have used the criterion of minimizing the probability of ruin of an insurance company to find the optimal rate of new business accumulation (Hipp and Taksar, 2000), to find the optimal investment strategy for the insurer (Hipp and Plum, 2000), and to find the optimal proportional reinsurance for the insurer (Schmidli, 2001).  Similarly, Olivieri and Pitacco (2003) consider the problem of maintaining solvency for a pension plan. 

From the individual investor's standpoint, which the view we take in this paper, Browne (1997) considers a financial model similar to ours (however, with no constraints on borrowing), and he maximizes the probability of reaching a safe level before ruining; also, see Browne (1995, 1999a, b) for related work.  This existence of the safe level is based on constant consumption and equals the price of a perpetuity that exactly covers the desired consumption.  He showed that no optimal policy exists for this problem because it is impossible to get to the safe region if one is maximizing this probability.  Our problem differs from his because we allow the individual to die before ruining, so that the individual can thereby ``win'' the game.  Plus, we optimize a different probability, namely the probability of ruining before dying.

In minimizing the probability that the individual outlives her money, we consider her random time of death, unlike Browne (1997). This assumption differs from the one usually assumed by financial planners and common retirement planning software in that they generally assume a specific age of death. A few earlier researchers have dealt with the problem of outliving one's wealth under the assumption of a random lifetime. For example, Milevsky, Ho, and Robinson (1997) and Milevsky and Robinson (2000) consider a random time of death modeled by using Canadian mortality data. They use simulation to find the probability of lifetime ruin. Our work differs from Milevsky, Ho, and Robinson (1997) and Milevsky and Robinson (2000) in that we find the optimal dynamic investment strategy to minimize the probability of lifetime ruin, while they take the investment strategy as fixed and find the corresponding probability of lifetime ruin.

Young (2004) considered the problem of minimizing the probability of lifetime ruin when the individual continuously consumes either a constant (real) dollar amount or a constant proportion of wealth.  For the most part, the results are intuitively appealing; as the model parameters vary, the changes in the ruin probability and the asset allocation are consistent with financial intuition.  However, when consumption is constant, she found that for wealth near 0, the optimal strategy is a heavily-leveraged position in the risky asset; that is, the individual borrows money to invest in the risky asset.  In order to avoid (nearly certain) ruin at low wealth levels, the individual takes on the lesser risk of borrowing a great deal of money at the riskless rate and investing it in the risky asset.  Although the objective of minimizing the probability of lifetime ruin is intuitively appealing, this leveraging at low wealth is not.

In this paper, we consider two model changes; one eliminates the leveraging, and the other reduces it.  In Section 2, we do not allow the individual to borrow money at all.  In Section 2.1, we present the financial model and define the corresponding probability of ruin.  In Section 2.2, we consider a constant rate of consumption, while in Section 2.3, we consider a rate of consumption that is proportional to wealth.  In the latter case, we show that our individual behaves exactly as an individual who maximizes expected discounted utility of consumption under a similar no-borrowing constraint and under power utility.  In Section 3, we parallel the work of Section 2 under the financial model that the individual can borrow money but only at a rate that is higher than the rate earned on the riskless asset.  In Section 4, we present a numerical example to demonstrate the results of Sections 2 and 3.  Section 5 concludes the paper.

\sect{2. Probability of Lifetime Ruin: No Borrowing}

In this section, we consider the problem of minimizing the probability of lifetime ruin when the individual is not allowed to borrow money.  We also impose the constraint that she cannot sell the risky asset short; however, this constraint is moot because she does not wish to short-sell the risky asset in the unconstrained case.  This occurs because we assume that the drift on the risky asset exceeds the rate of return on the riskless asset.  In Section 2.1, we present the financial market and define the probability of ruin.  In Section 2.2, we solve for the probability of ruin in the case for which consumption is constant.  Section 2.3 parallels Section 2.2 in the case for which consumption is proportional to wealth.  We show how this latter case is related to optimal investment and consumption in a Merton model with no borrowing.

\subsect{2.1. Financial Market}

In this section, we first present the financial ingredients that make up the individual's wealth, namely, consumption, a riskless asset, and a risky asset. We, then, define the minimum probability of lifetime ruin under the constraint of no borrowing (and no short-selling, although this is not a real restriction).

The individual consumes at a continuous rate $c(w)$, in which $w$ is her current wealth. In this paper, we consider two forms of the consumption function:

\item{(1)} $c(w) = c$; the individual consumes at a constant (nominal or real) dollar rate $c$. If $c$ is a real dollar rate, then returns in the financial market are real, too.

\item{(2)} $c(w) = pw$; the individual consumes a constant proportion $p$ of her wealth.

\smallskip

We assume that the individual invests in a riskless asset whose price at time $t$, $X_t$, follows the process $dX_t = rX_t dt, X_0 = x > 0$, for some fixed rate of interest $r > 0$. Also, the individual invests in a risky asset whose price at time $t$, $S_t$, follows geometric Brownian motion given by

$$\left\{ \eqalign{dS_t &= \mu S_t dt + \sigma S_t dB_t, \cr
S_0 &= S > 0,} \right. \eqno(2.1)$$

\noindent in which $\mu > r$, $\sigma > 0$, and $B$ is a standard Brownian motion with respect to a filtration of the probability space $(\Omega, {\cal F}, {\bf Pr})$.

Let $W_t$ be the wealth at time $t$ of the individual, and let $\pi_{0,t}$ be the amount that the decision maker invests in the risky asset at that time.  We use a subscript $0$ to denote the fact that no borrowing (or short-selling) is allowed.  It follows that the amount invested in the riskless asset is $W_t - \pi_{0,t}$. Thus, wealth follows the process

$$\left\{ \eqalign{dW_t &= [rW_t + (\mu - r) \pi_{0,t}  - c(W_t)] dt + \sigma \pi_{0,t} dB_t, \cr
W_0 &= w.} \right. \eqno(2.2)$$
			
By ``outliving her wealth,'' or equivalently ``lifetime ruin,'' we mean that the individual's wealth reaches zero before she dies.  Let $\tau_0$ denote the first time that wealth equals zero, and let $\tau_d$ denote the random time of death of our individual.  We assume that $\tau_d$ is exponentially distributed with parameter $\lambda$ (that is, with expected time of death equal to $1/\lambda$); this parameter is also known as the {\it hazard rate} of the individual.
			
Denote the minimum probability that the individual outlives her wealth by $\psi_0(w)$, in which the argument $w$ indicates that one conditions on the individual possessing wealth $w$ at the current time.  Recall that the subscript $0$ reminds us that no borrowing is allowed.  Thus, $\psi_0$ is the minimum probability that $\tau_0 < \tau_d$, in which one minimizes with respect to admissible investment strategies $\pi_0$.  A strategy $\pi_0$ is {\it admissible} if it is ${\cal F}_t$-progressively measurable (in which ${\cal F}_t$ is the augmentation of $\sigma(W_s: 0 \le s \le t)$) and if it satisfies the integrability condition $\int_0^t \pi_{0,s}^2 ds < \infty$, almost surely, for all $t \ge 0$.  Also, in this section, we restrict $\pi_{0,t} \in [0, W_t]$; that is, no borrowing nor short-selling is allowed.  Because $\mu > r$, the latter is effectively not a restriction.  It follows that one can express $\psi_0$ as follows:

$$\psi_0(w) = \inf_{\pi_0} {\bf Pr} \left[\tau_0 < \tau_d \vert W_0 = w \right]. \eqno(2.3)$$

For every $\alpha \in \bf R$, we associate a second-order differential operator ${\cal L}^{\alpha}$ with this minimization problem as follows:  For every open set $G \subset {\bf R}^+$ and for every $h \in C^2(G)$, we define the function ${\cal L}^{\alpha} h: G \rightarrow \bf R$ by

$${\cal L}^{\alpha} h(w) = \left[ rw + (\mu - r) \alpha - c(w) \right] h'(w) + {1 \over 2} \sigma^2 \alpha^2 h''(w) - \lambda h(w). \eqno(2.4)$$

\noindent We use ${\cal L}^{\alpha}$ in the following sections to characterize the minimum probability of ruin $\psi_0$ in the cases for which $c(w) = c$ and $c(w) = pw$, respectively.

\subsect{2.2. Constant Consumption}

In this section, we consider the case for which the consumption rate $c(w) = c$, a positive constant.  Note that when wealth reaches $c/r$, then the individual can place all her wealth in the riskless asset and consume $c$ continuously without risk of running out of money.  Therefore, the probability of lifetime ruin equals 0 when wealth is greater than or equal to $c/r$.

Define the stopping time $\tau = \tau_d \wedge \tau_{c/r}$, in which $\tau_{c/r} = \inf\{t > 0: W_t = c/r \}$, with the convention that $\inf \emptyset = \infty$.  It follows that $\psi_0(w) = \inf_{\pi_0 \in {\cal A}} {\bf Pr} \left[ \tau_0 < \tau \vert W_0 = w \right]$, in which ${\cal A}$ is the set of admissible strategies.  We have the following verification theorem.

\th{2.1} {Suppose $h_0$ is a decreasing function from ${\bf R}^+$ to $[0, 1]$ and suppose $\alpha_0$ is a function from ${\bf R}^+$ to ${\bf R}^+$ that satisfy the following conditions:
\item{$(i)$} $h_0 \in C^2$ on $[0,c/r);$
\item{$(ii)$} $\alpha_0 \in {\cal A};$
\item{$(iii)$} ${\cal L}^{\alpha} h_0(w) \ge 0,$ for $0 \le \alpha \le w < c/r;$
\item{$(iv)$} ${\cal L}^{\alpha_0(w)} h_0(w) = 0,$ for $w \in (0, c/r);$
\item{$(v)$} $h_0(0) = 1,$ and $h_0(w) = 0$ for $w \ge c/r$.}

{\it Under the above conditions, the minimum probability of the lifetime ruin $\psi_0$ is given by}
$$\psi_0(w) = h_0(w), \quad w \ge 0, \eqno(2.5)$$

\noindent {\it and the optimal investment strategy in the risky asset $\pi^*_0$ is given by}
$$\pi^*_0(w) = \alpha_0(w), \quad w \in [0, c/r]. \eqno(2.6)$$

\medskip

\pf  Assume that we have $h_0$ as specified in the statement of this theorem.  Let $N$ denote a Poisson process with rate $\lambda$ such that $N$ independent of the Brownian motion $B$ driving the wealth.  The first occurrence of $N$ represents the death of the individual.

Let $\alpha$ be some function on $[0, c/r]$ such that $0 \le \alpha(w) \le w$.  Let $W^\alpha$ denote the wealth process when we use $\alpha$ as the investment policy, and let $\alpha_s = \alpha(W^\alpha_s)$.  Let $\Delta$ represent the ``coffin state," and let ${\bf R}^+ \cup \{\Delta\} $ be the one point compactification of ${\bf R}^+$.  For any function $f:{\bf R}^+ \rightarrow {\bf R}^+$, we define its extension to $ {\bf R}^+ \cup \Delta$ by assigning $f(\Delta) = 0$. We kill the wealth process as soon as the Poisson process jumps (that is, when the individual dies) and assign $W_{\tau_d} = \Delta$.   Note that $h_0(c/r) = h_0(W^\alpha_{\tau_d}) = 0$.  

By applying It\^o's lemma, we have

$$\eqalign{& h_0(W^\alpha_{t \wedge \tau_0 \wedge \tau}) =h_0(w)+ \int_{0}^{t \wedge \tau_0 \wedge \tau} \left( (rW^\alpha_s + (\mu-r) \alpha_s - c)h_0'(W^\alpha_s) + {1 \over 2} \sigma^2 \alpha^2_s h_0''(W^\alpha_s) \right) ds \cr
& \qquad \qquad \qquad - \lambda \int_0^{t \wedge \tau_0 \wedge \tau} h_0(W^\alpha_{s-})ds+ \int_0^{t \wedge \tau_0 \wedge \tau} h_0'(W^\alpha_s)\sigma \alpha_s dB_s \cr
& \qquad \qquad \qquad + \int_0^{t \wedge \tau_0 \wedge \tau} (h_0(W^\alpha_s)-h_0(W^\alpha_{s-})) d(N_s-\lambda s) \cr
& = h_0(w) + \int_0^{t \wedge \tau_0 \wedge \tau}{\cal L}^{\alpha}h_{0}(W^\alpha_s) ds + \int_{0}^{t \wedge \tau_0 \wedge \tau} h_0'(W^\alpha_s)\sigma \alpha_s dB_{s} -h_0(W^\alpha_{s-}) d(N_s-\lambda s).} \eqno(2.7)$$

\noindent The second equality follows from the definition of ${\cal L}^{\alpha}$ and the fact that $h_0(W^\alpha_{\tau_d}) = 0$. Now if we take the expectation of both sides, the expectation of the third term in (2.7) becomes zero because

$$\eqalign{& {\bf E}_{w} \left[ \int_0^{t \wedge \tau_{0} \wedge \tau} \left((h_{0}'(W^{\alpha}_{s}))^{2} \sigma^{2} \alpha_{s}^{2}+ \lambda h_{0}(W^{\alpha}_{s})\right)ds \right] \cr
& \qquad < {\bf E}_{w} [ \tau_d] (\max_{0 \leq w \leq c/r} (h'_{0}(w))^{2} \sigma^{2} (c^{2}/r^{2})+\lambda) < \infty,} \eqno(2.8)$$

\noindent since ${\bf E}_w[\tau_d] = 1/\lambda$, and $h'(w)$ is bounded on $[0, c/r)$ by assumption (i).  The first inequality in (2.8) follows from assumption (ii) of the theorem.  ${\bf E}_w$ denotes the conditional expectation given $W_0=w$.

Then, we have

$${\bf E}_{w} [h_0 (W^{\alpha}_{t \wedge \tau_0 \wedge \tau})]= h_0(w)+ {\bf E}_w\left[\int_0^{t \wedge \tau_0 \wedge
\tau} {\cal L}^{\alpha}h_{0}(W^{\alpha}_s) ds \right] \ge h_{0}(w), \eqno(2.9)$$

\noindent where the inequality follows from assumption (iii) of the theorem.  The expression in (2.9) shows that $\left( h_0(W^{\alpha}_{t \wedge \tau_0 \wedge \tau}) \right)_{t \ge 0}$ is a sub-martingale.

Now, since $h_{0}(0) = 1$, $h_{0}(W^\alpha_{\tau_d}) = 0$, and $h_{0}(W^{\alpha}_{c/r})=0$, we have

$$h_{0}(W^{\alpha}_{\tau_0 \wedge \tau}) = {\bf 1}_{\{ \tau^\alpha_0< \tau\}}, \eqno(2.10)$$

\noindent in which {\bf 1} is the indicator function. By taking expectations of both sides of (2.10), we have

$${\bf E}_w[h_0(W^{\alpha}_{\tau_0 \wedge \tau})] = {\bf Pr}_w(\tau^\alpha_0 < \tau) \ge h_0(w). \eqno(2.11)$$

\noindent ${\bf Pr}_w$ denotes the conditional probability given $W_0 = w$. We write $\tau^\alpha_0$ to emphasize the dependence of $\tau_0$ on the strategy $\alpha$.  The last inequality follows from an aplication of optional sampling theorem since $\left( h_0(W^\alpha_{t \wedge \tau_0 \wedge \tau}) \right)_{t \ge 0}$ is a sub-martingale.  (We can apply optional sampling theorem due to Theorem 3.15 of Karatzas and Shereve (1991) since $h_0(w) \in [0,1]$ for all $w \ge 0$.)  Therefore,

$$\inf_{\alpha \in {\cal A}} {\bf Pr}_w(\tau^\alpha_0< \tau) = \psi(w) \ge h_{0}(w). \eqno(2.12)$$

Let $\alpha_0$ be as specified in the statement of this theorem; that is, $\alpha_0$ is the minimizer of ${\cal L}^\alpha h_0$.   It follows that $\left(h_0(W^{\alpha_0}_{t \wedge \tau_0 \wedge \tau}) \right)_{t \ge 0}$ is a martingale.  Therefore,

$${\bf E}_w [h_0 (W^{\alpha_0}_{\tau_0 \wedge \tau})] = {\bf Pr}_w( \tau_0^{\alpha_0}< \tau) = h_0(w); \eqno(2.13)$$

\noindent hence, we have demonstrated (2.5) and (2.6) on $[0, c/r)$.  Assumption (v) and the remark immediately preceding the statement of this theorem complete the proof.  \hfill {\bf QED}

\medskip

At this point, we could simply write down $h_0$ and $\alpha_0$ and prove a theorem that they satisfy the hypotheses of Theorem 2.1, and we would be done.  However, that would not show the reader how to solve such optimization problems.  Therefore, in the remainder of this section, we demonstrate the steps of how to find $h_0$ and $\alpha_0$, and we summarize our results at the end.

In the unconstrained case (Young, 2004), the minimum probability of ruin $\psi$ is given by

$$\psi(w) = \left( 1 - {r \over c} w \right)^d,  \eqno(2.14)$$

\noindent with

$$d = {1 \over 2r} \left[ (r + \lambda + m) + \sqrt{(r + \lambda + m)^2 - 4r \lambda} \right] > 1, \eqno(2.15)$$ and

$$m = {1 \over 2} \left( {\mu - r \over \sigma} \right)^2. \eqno(2.16)$$

Because $\psi$ is convex, the corresponding unconstrained optimal investment in the risky asset $\pi^*$ is given by the first-order necessary condition implicit in items (iii) and (iv) of Theorem 2.1:

$$\pi^*(w) = {\mu - r \over \sigma^2} {1 \over d - 1} \left( {c \over r} - w \right), \eqno(2.17)$$

\noindent a non-negative, decreasing, linear function of wealth.  Note that as wealth increases towards $c/r$, the amount invested in the risky asset decreases to zero. This makes sense because as the individual becomes wealthier, she does not need to take on as much risk to achieve her fixed consumption rate of $c$.

On the other hand, for wealth small, the (unconstrained) optimal amount invested in the risky asset is greater than wealth; that is, the individual borrows money to invest in the risky asset in order to avoid the greater risk of lifetime ruin.  We believe that most people with small wealth will not borrow to invest in a risky asset to avoid ruin and that no credible financial advisor will give such advice.  Therefore, it is reasonable to include the restriction that the amount invested in the risky asset cannot exceed current wealth.

In order to find $h_0$ and $\alpha_0$ that satisfy Theorem 2.1, we hypothesize that they satisfy additional properties not explicitly stated in Theorem 2.1.  If we find $h_0$ and $\alpha_0$ that satisfy Theorem 2.1 and the additional properties, then these additional properties are implicit in Theorem 2.1 because of the uniqueness of $h_0$ on $(0, c/r]$.  The function $h_0$ is unique, if it exists, because Theorem 2.1 states that if the Hamilton-Jacobi-Bellman equation (defined by items (i) through (v)) has a solution, then that solution equals the minimum probability $\psi_0$.  We propose the following ansatz.

{\bf Ansatz:} By considering the amount invested in the {\it riskless} asset in the unconstrained case, namely $w - \pi^*(w)$, we conjecture the form of the constrained $\alpha_0$.  To that end, note that $w - \pi^*(w)$ is an increasing, linear function of wealth; also, the amount invested in the riskless asset for $w = 0$ is negative and for $w = c/r$ is positive.  Therefore, we hypothesize that the optimal {\it constrained} amount invested in the riskless asset is (1) a continuous function of wealth; (2) 0 for wealth below some level, say $w_l$, and (3) positive for wealth greater than $w_l$.  The subscript $l$ on $w_l$ stands for {\bf l}ending because for wealth above $w_l$, the individual lends money to the financial institution selling the riskless asset.  In other words, we suppose that the optimal investment in the riskless asset under the constraint of no borrowing is a truncated version of the optimal investment in the unconstrained case.  We are not asserting that the wealth level above which the individual invests a positive amount in the riskless asset $w_l$ is the same in both cases; however, that turns out to be the case, as we note below.  In addition to these properties of $\alpha_0$, we hypothesize that $h_0$ is {\it strictly} decreasing on $[0, c/r]$. 

\medskip

Consider each of the two intervals $(w_l, c/r]$ and $[0, w_l]$ in turn.  First, consider $(w_l, c/r]$, on which we hypothesize that the borrowing constraint does not bind.  We have the following proposition concerning $h_0$ on this interval.

\prop{2.2} {Under the hypothesis that $0 \le \alpha_0(w) < w$ on $(w_l, c/r],$ the function $h_0$ given by}

$$h_0(w) = \beta_0 \left( 1 - {r \over c} w \right)^d, \eqno(2.18)$$

\noindent{\it with $d$ as in $(2.15)$ and $\beta_0 \ge 1$, satisfies the properties specified in Theorem 2.1.  Thus, on $(w_l, c/r]$, $h_0$ is a multiple of the unconstrained probability of ruin $\psi$.  The corresponding $\alpha_0$ on this interval is given by}

$$\alpha_0(w) = {\mu - r \over \sigma^2} {1 \over d - 1} \left( {c \over r} - w \right),  \eqno(2.19)$$

\noindent {\it identical to $\pi^*$ in the unconstrained case.}

\medskip

\pf From items (iii), (iv), and (v) of Theorem 2.1, we solve

$$\lambda h_0(w) = (rw - c) h'_0(w) + \min_{\alpha} \left[ (\mu - r) \alpha h'_0(w) + {1 \over 2} \sigma^2 \alpha^2 h''_0(w) \right], \eqno(2.20)$$

\noindent with boundary conditions $h_0(c/r) = 0$.  We also have the boundary condition $h'_0(c/r)$, which we demonstrate now.  Consider the solution $\phi$ of (2.20) with $\lambda = 0$.  Certainly, we have $h_0 \le \phi$ on an interval $(c/r - \delta, c/r]$ for some $\delta > 0$ because the probability of ruining before dying is no greater than the probability of ruining before infinity.  If we show that $\phi'(c/r) = 0$, then $h'_0(c/r)$ because $h_0$ is wedged between 0 and $\phi$ as wealth approaches 0.  Note that $\phi$ solves the equation

$$0 = (rw - c) \phi'(w) + \min_{\alpha} \left[ (\mu - r) \alpha \phi'(w) + {1 \over 2} \sigma^2 \alpha^2 \phi''(w) \right], \eqno(2.21)$$

\noindent with $\phi(c/r) = 0$.  From Pestien and Sudderth (1985), the optimal investment strategy $\alpha^*$ is the one that maximizes the ratio of the drift to the variance, or equivalently,

$$f(\alpha) = {(\mu - r)\alpha - (c - rw) \over \alpha^2}. \eqno(2.22)$$

\noindent The value of $\alpha$ that maximizes this expression is $\alpha^*(w) = 2(c - rw)/(\mu - r)$.  On the other hand, we also have that

$$\alpha^*(w) = -{\mu - r \over \sigma^2} {\phi'(w) \over \phi''(w)} \eqno(2.23)$$

\noindent from the first-order condition necessary in (2.21).  Thus,

$$\phi'(w) = k (c - rw)^{m/r}, \quad w \in (c/r - \delta, c/r], \eqno(2.24)$$

\noindent for some $k < 0$.  Thus, $\phi'(c/r) = 0$, from which it follows that $h'_0(c/r) = 0$.

Now, if $h_0$ is not convex in some neighborhood of a point $w^* \in (w_l, c/r]$, then $h''_0(w) < 0$ on that neighborhood of $w^*$.  It follows that the optimal solution $\alpha_0$ is as large as possible on that neighborhood, contradicting our hypothesis that the borrowing constraint does not hold on $(w_l, c/r]$.  Therefore, to be consistent with this hypothesis, $h_0$ is convex on $(w_l, c/r]$, and we consider the Legendre transform $\tilde h$ of $h_0$ defined by

$$\tilde h(v) = \min_w [h_0(w) + wv]. \eqno(2.25)$$

\noindent Note that we can recover $h_0$ from $\tilde h$ by

$$h_0(w) = \max_v [\tilde h(v) - wv]. \eqno(2.26)$$

\noindent The minimizing value of $w$ in (2.25) equals $I(-v) = \tilde h'(v)$, in which $I$ is the inverse function of $h'_0$.  Therefore, the maximizing value of $v$ in (2.26) equals $-h'_0(w)$.

Substitute $w = I(-v)$ in equation (2.20) to obtain

$$\lambda \tilde h(v) + (r - \lambda) v \tilde h'(v) - m v^2 \tilde h''(v) = cv, \eqno(2.27)$$

\noindent in which $m$ is given in (2.16).  The general solution of (2.27) is

$$\tilde h(v) = \tilde D_1 v^{\tilde B_1} + \tilde D_2 v^{\tilde B_2} + {c \over r} v, \eqno(2.28)$$

\noindent in which $\tilde D_1$ and $\tilde D_2$ are constants to be determined, and $\tilde B_1$ and $\tilde B_2$ are the positive and negative roots, respectively of

$$-\lambda - (r - \lambda + m) \tilde B + m \tilde B^2 = 0.  \eqno(2.29)$$

\noindent Thus,

$$\tilde B_1 = {1 \over 2 m} \left[ (r - \lambda + m) + \sqrt{(r - \lambda + m)^2 + 4 \lambda m} \right] > 1, \eqno(2.30)$$

\noindent and

$$\tilde B_2 = {1 \over 2 m} \left[ (r - \lambda + m) - \sqrt{(r - \lambda + m)^2 + 4 \lambda m} \right] < 0. \eqno(2.31)$$

Define $v_c = -h'_0(c/r) = 0$; that is, $\tilde h'(0) = c/r.$  From the definition of $\tilde h$ in (2.25) and from $h_0(c/r) = 0$, at $v = v_c = 0$, we have

$$\tilde h(0) = 0. \eqno(2.32)$$

\noindent It follows that $\tilde D_2 = 0$.  We can, then, recover $h_0$ from (2.28) and (2.26) and obtain the expression for $h_0$ in (2.18).

Because $h_0$ is convex on $(w_l, c/r]$, the optimal policy $\alpha_0$ is given by the first-order necessary condition in (2.20).  This observation leads to the expression for $\alpha_0$ given in equation (2.19).  \hfill {\bf QED}

\medskip

\cor{2.3} {By the assumed continuity of $\alpha_0,$}

$$w_l = {x \over 1 + x} {c \over r}, \eqno(2.33)$$

\noindent {\it in which}

$$x = {\mu - r \over \sigma^2} {1 \over d - 1}.  \eqno(2.34)$$

\medskip

Because $\alpha_0 = \pi^*$ on $[w_l, c/r]$, the lending level $w_l$ is the same in both the unconstrained case and in the no-borrowing case, in which the lending level is the wealth level above which the individual invests a positive amount in the riskless asset.  We find this result rather surprising because we expected that the individual would have a higher lending level in the constrained case in order to invest more money in the risky asset due to the fact that she cannot invest as much in the risky asset as she wishes when her wealth is below $w_l$.  We will see this myopia again in Section 3 when we consider the case of borrowing at a rate higher than the lending rate.

Next, consider the interval $[0, w_l]$, on which $\alpha_0(w) = w$.  We have the following proposition concerning $h_0$ on this interval.  We state it without proof because it easily follows from Theorem 2.1 and from our hypothesis that $\alpha_0(w) = w$ on $[0, w_l]$.

\prop{2.4} {Under the hypothesis that $\alpha_0(w) = w$ on $[0, w_l]$, the function $h_0$ solves}

$$\lambda h_0 = (\mu w - c) h_0' + {1 \over 2} \sigma^2 w^2 h_0'',  \eqno(2.35)$$

\noindent {\it with boundary conditions}

$$h_0(0) = 1 \eqno(2.36)$$

\noindent {\it and}

$${h_0(w_l) \over h'_0(w_l)} = -{1 \over d} \left( {c \over r} - w_l \right). \eqno(2.37)$$

\medskip

These two boundary conditions allow us to solve the second-order ode (2.35) numerically.  Once we have a solution for $h_0$ on $[0, w_l]$, then we can use the condition $h_0(w_l-) = h_0(w_l+)$ to determine the remaining unknown parameter $\beta_0$.

It remains for us to show that if $h_0$ is given as the solution of $(2.35) - (2.37)$, then the constrained optimal investment in the risky asset is equal to the current wealth.  We state this formally in the following proposition.

\prop{2.5} {Suppose $h_0$ on $[0, w_l]$ solves $(2.35) - (2.37)$, then}

$$\arg \min_{0 \le \alpha \le w} \left[ (\mu - r) \alpha h'_0(w) + {1 \over 2} \sigma^2 \alpha^2 h''_0(w) \right] = w, \quad 0 \le w \le w_l. \eqno(2.38)$$

\pf  Define the parabola $f$ by

$$f(\alpha) = (\mu - r) \alpha h'_0(w) + {1 \over 2} \sigma^2 \alpha ^2 h''_0(w), \eqno(2.39)$$

\noindent for $w \in [0, w_l]$.  Proving this proposition is equivalent to showing that

$$f'(w) = (\mu - r) h'_0(w) + \sigma^2 w h''_0(w) \le 0, \eqno(2.40)$$

\noindent for all wealth in $[0, w_l]$.  By substituting for $h''_0$ from the differential equation (2.35), this inequality becomes

$$[-(\mu + r) w + 2c] h'_0(w) + 2 \lambda h_0(w) \le 0. \eqno(2.41)$$

Towards the goal of demonstrating inequality (2.41) on $[0, w_l]$, define the function $y$ by

$$y(w) = {h_0(w) \over h'_0(w)}, \eqno(2.42)$$

\noindent for $w \in [0, w_l]$.  Note that showing inequality (2.41) is equivalent to showing that $y \ge z$ on $[0, w_l]$, in which $z$ is given by

$$z(w) = {\mu + r \over 2 \lambda} w - {c \over \lambda}. \eqno(2.43)$$

\noindent We complete the proof of this proposition by proving the following lemma. 

\lem{2.6} {$y > z$ on $(0, w_l),$ with equality at $w = 0$ and $w = w_l$.}

\medskip

\pf By substituting in the linear second-order differential equation (2.35), we obtain a non-linear first-order differential equation for $y$ on $[0, w_l]$:

$$\sigma^2 w^2 (y'(w) - 1) = -2 \lambda y^2(w) + 2(\mu w - c) y(w). \eqno(2.44)$$

\noindent The function $z$ given by (2.43) solves a similar differential equation:

$$\sigma^2 w^2 \left( z'(w) - { \mu + r \over 2 \lambda} \right) = -2 \lambda z^2(w) + 2 \left( { \mu + r \over 2} w - c \right) z(w). \eqno(2.45)$$

Note that $y(0) = z(0) = -c/\lambda$ and $y(w_l) = z(w_l) = -(1/d) (c/r - w_l)$.  We first show that $y'(w_l) < z'(w_l)$.  Indeed, after substituting for $y(w_l)$ in (2.44) and simplifying, we obtain

$$y'(w_l) = 1 + {r + m \over \lambda} - {r \over \lambda}d. \eqno(2.46)$$

\noindent By substituting for $d$ in (2.46), we have that $y'(w_l) < z'(w_l)$ if and only if

$$-\mu + \lambda + m < \sqrt{(r + \lambda + m)^2 - 4 r \lambda}, \eqno(2.47)$$

\noindent which, because $\mu > r$, holds if

$$-r + \lambda + m < \sqrt{(r + \lambda + m)^2 - 4 r \lambda}. \eqno(2.48)$$

\noindent Inequality (2.48) is true as one can see by squaring both sides and simplifying.  Thus, we have shown that $y'(w_l) < z'(w_l)$.

It follows that $y > z$ on $(w_l - \delta, w_l)$ for some $\delta > 0$. Next, suppose that there exists $w^* \in (0, w_l)$ such that $y(w^*) = z(w^*)$ and $y > z$ on $(w^*, w_l)$.  We will show that such a $w^*$ cannot exist.  

Because $y(w^*) = z(w^*)$ and $y > z$ on $(w^*, w_l)$, we have $y'(w^*) \ge z'(w^*)$, which implies that

$$\eqalign{& 1 - {2 \lambda \over \sigma^2 (w^*)^2} y(w^*)^2 + {2(\mu w^* - c) \over \sigma^2 (w^*)^2} y(w^*) \cr
& \quad \ge {\mu + r \over 2 \lambda} - {2 \lambda \over \sigma^2 (w^*)^2} z(w^*)^2 + {2 \left( {\mu + r \over 2} w^* - c \right) \over \sigma^2 (w^*)^2} z(w^*) \cr
\Rightarrow & 1 + {2(\mu w^* - c) \over \sigma^2 (w^*)^2} z(w^*) \ge {\mu + r \over 2 \lambda} + {2 \left( {\mu + r \over 2} w^* - c \right) \over \sigma^2 (w^*)^2} z(w^*) \cr
\Rightarrow & 1 - {\mu + r \over 2 \lambda} \ge - { \mu - r \over \sigma^2 w^*} \left( {\mu + r \over 2 \lambda} w^* - {c \over \lambda} \right).} \eqno(2.49)$$

\noindent Note that the right-hand-side of the last line of (2.49) is positive.  Therefore, if ${\mu + r \over 2 \lambda} \ge 1$, then we have our contradiction.

To continue, suppose ${\mu + r \over 2 \lambda} < 1$.  In that case, the inequality in (2.49) becomes

$$w^* \ge {2c(\mu - r) \over \sigma^2(2 \lambda - \mu - r) + (\mu^2 - r^2)}. \eqno(2.50)$$

\noindent If $w_l < {2c(\mu - r) \over \sigma^2(2 \lambda - \mu - r) + (\mu^2 - r^2)}$, then $w^* \in (0, w_l)$ cannot exist, and we are done.  Now, by substituting for $w_l$ and simplifying, one can show that $w_l < {2c(\mu - r) \over \sigma^2(2 \lambda - \mu - r) + (\mu^2 - r^2)}$ is equivalent to inequality (2.47), which we have already demonstrated.  Therefore, no such $w^*$ exists, and $y > z$ on $(0, w_l)$.  \hfill {\bf QED}

\medskip

We have the following theorem that summarizes the results of this section.

\th{2.7} {For constant consumption $c(w) = c$, the constrained minimum probability of ruin $\psi_0 \in C^{1}({\bf R^{+}})  \cap C^2( {\bf R}^+ - \{ c/r \}) $ is given by}

$$\psi_0(w) = h_0(w), \quad 0 \le w \le w_l, \eqno(2.51)$$

\noindent {\it in which $h_0$ solves $(2.35) - (2.37)$ with $w_l$ in $(2.33),$ and}

$$\psi_0(w) = \beta_0 \left( 1 - {r \over c} w \right)^d, \quad w_l < w \le c/r, \eqno(2.52)$$

\noindent {\it in which $\beta_0 = h_0(w_l)(1 - rw_l/c)^{-d}$.  Finally, $\psi_0(w) = 0$ for $w > c/r$.

The optimal investment strategy is given by}

$$\pi^*_0(w) = \cases{w, & if $0 \le w \le w_l;$ \cr {\mu - r \over \sigma^2} {1 \over d - 1} \left( {c \over r} - w \right), & if $w_l < w \le c/r$.} \eqno(2.53)$$

\medskip

One can show that if initial wealth lies below $c/r$, then it never reaches that ``safe'' level; that is, $\tau^{\alpha_0}_{c/r} = \infty$.  See Young (2004) and Browne (1997) for discussions of this phenomenon in related models.  In Section 4, we present numerical examples that demonstrate the results of this section.

\subsect{2.3. Proportional Consumption}

In this section, we consider the case for which the consumption rate is proportional to wealth $c(w) = pw$, in which $p > r$.  Note that if we specify that ruin occurs if wealth reaches 0, then the individual with this consumption function does not ruin.  Therefore, for this case, we say that ruin occurs when wealth reaches $w_0 > 0$; that is, let $\tau_0$ in definition (2.3) be the (first) time that wealth reaches $w_0$.

We have the following verification theorem whose proof we omit because it is similar to the proof for Theorem 2.1.

\th{2.8} {Suppose $g_0$ is a strictly decreasing function from ${\bf R}^+$ to $[0, 1]$ and suppose $\gamma_0$ is a function from ${\bf R}^+$ to ${\bf R}$ that satisfy the following conditions:
\item{$(i)$} $g_0 \in C^2 \left( {\bf R}^+ \right);$
\item{$(ii)$} $0 \le \gamma_0(w) \le w$ for $w \ge 0;$
\item{$(iii)$} ${\cal L}^{\gamma} g_0(w) \ge 0,$ for $0 \le \gamma \le w;$
\item{$(iv)$} ${\cal L}^{\gamma_0(w)} g_0(w) = 0,$ for $w \ge 0;$
\item{$(v)$} $g_0(0) = 1,$ and $\lim_{w \rightarrow \infty} g_0(w) = 0$.}

{\it Under the above conditions, the minimum probability of the lifetime ruin $\psi_0$ is given by}
$$\psi_0(w) = g_0(w), \quad w \ge 0, \eqno(2.54)$$

\noindent {\it and the optimal investment strategy in the risky asset $\pi^*_0$ is given by}
$$\pi^*_0(w) = \gamma_0(w), \quad w \ge 0. \eqno(2.55)$$

\medskip

In the unconstrained case (Young, 2004), the minimum probability of ruin $\psi$ is given by

$$\psi(w) = \left( {w \over w_0} \right)^{-a_r},  \eqno(2.56)$$

\noindent with

$$a_r = {1 \over 2(p - r)} \left[ (r - p + \lambda + m) + \sqrt{(r - p + \lambda + m)^2 + 4 \lambda (p - r)} \right] > 0, \eqno(2.57)$$

\noindent and $m$ is as in (2.16).  Note that $-a_r$ equals $d$ in (2.15) with $r$ replaced by $r - p$.

Because $\psi$ is convex, the corresponding unconstrained optimal investment in the risky asset $\pi^*$ is given by the first-order necessary condition implicit in items (iii) and (iv) of Theorem 2.8:

$$\pi^*(w) = {\mu - r \over \sigma^2} {w \over a_r + 1}, \eqno(2.58)$$

\noindent a positive proportion of wealth.  Note that $\pi^*(w) < w$ if and only ${\mu - r \over \sigma^2} {1 \over a_r + 1} < 1$.  This observation corresponds neatly with the following theorem.

\th{2.9} {The constrained minimum probability of ruin is given by}

$$\psi_0(w) = \cases{ \left( {w \over w_0} \right)^{-a_r}, & if ${\mu - r \over \sigma^2} {1 \over a_r + 1} < 1;$ \cr
\left( {w \over w_0} \right)^{-k}, & if ${\mu - r \over \sigma^2} {1 \over a_r + 1} \ge 1,$} \eqno(2.59)$$

\noindent {\it in which}

$$k = {1 \over \sigma^2} \left[ \left( \mu - p - {1 \over 2} \sigma^2 \right) + \sqrt{\left( \mu - p - {1 \over 2} \sigma^2 \right)^2 + 2 \sigma^2 \lambda} \right] > 0. \eqno(2.60)$$

\noindent {\it The corresponding constrained optimal investment strategy is given by}

$$\pi^*_0(w) = \cases{ {\mu - r \over \sigma^2} {w \over a_r + 1}, & if ${\mu - r \over \sigma^2} {1 \over a_r + 1} < 1;$ \cr
w, & if ${\mu - r \over \sigma^2} {1 \over a_r + 1} \ge 1.$} \eqno(2.61)$$

\medskip

\pf Note that if $g(w)$ solves the equation

$$\lambda g(w) = (p - r) w g'(w) + \min_{0 \le \gamma \le w} \left[ (\mu - r) \gamma g'(w) + {1 \over 2} \sigma^2 \gamma^2 g''(w) \right], \eqno(2.62)$$

\noindent then $g(\alpha w)$ and $\beta g(w)$ also solve the equation for any $\alpha > 0$ and $\beta > 0$.  Because $g_0(w_0) = 1$, we have that $g_\alpha$ defined by $g_\alpha(w) = g_0(\alpha w)$ solves (2.62) with $g_\alpha(w_0/\alpha) = 1$.  Also, $\hat g_\beta$ defined by $\hat g_\beta(w) = \beta g_0(w)$ solves (2.62) with $\hat g_\beta(w_0/\alpha) = \beta g_0(w_0/\alpha)$.  Set $\beta = 1/g_0(w_0/\alpha)$; then, $\hat g_\beta(w_0/\alpha) = 1$, and by uniqueness of $g_0$, we have that $g_\alpha(w) = \hat g_\beta(w)$ for all $w > 0$.  We, therefore, have demonstrated the following functional equation for $g_0$:

$$g_0(uv) = g_0(u) g_0(v w_0), \quad u, v > 0. \eqno(2.63)$$

Define $\rho(w) = g_0(w w_0)$; then, $\rho$ solves

$$\rho(uv) = \rho(u) \rho(v), \quad u, v > 0. \eqno(2.64)$$

\noindent Under mild regularity conditions on $\rho$, such as left- or right-continuity, it is well-known that $\rho$ is a power function, say $w^{-a}$, for some $a \in {\bf R}$.  Thus, $g_0$ is given by

$$g_0(w) = \left( {w \over w_0} \right)^{-a}, \eqno(2.65)$$

\noindent for $a > 0$ because $g_0$ is decreasing.

All that remains is for us to determine the value of $a$ under various ranges of parameters.  Items (iii) and (iv) in Theorem 2.8 lead us to solve

$$\eqalign{\lambda \left( {w \over w_0} \right)^{-a} &= (p - r) a \left( {w \over w_0} \right)^{-a} \cr
&+ \min_{0 \le \gamma \le w} \left[ (\mu - r) \gamma \left( {-a \over w_0} \right) \left( {w \over w_0} \right)^{-a-1} + {1 \over 2} \sigma^2 \gamma^2 {a(a + 1) \over w_0^2} \left( {w \over w_0} \right)^{-a-2}  \right],} \eqno(2.66)$$

\noindent which reduces to

$$\lambda = a (p - r) + a \min_{0 \le \gamma \le w} \left[ - (\mu - r) {\gamma \over w} + {1 \over 2} \sigma^2 (a + 1) {\gamma^2 \over w^2}  \right]. \eqno(2.67)$$

\noindent The quantity in square brackets is minimized by $\gamma = w$ if ${\mu - r \over \sigma^2} {1 \over a_r + 1} \ge 1$; otherwise, it is minimized by $\gamma = {\mu - r \over \sigma^2} {w \over a_r + 1}$.  In the former case, $a = a_r$ in (2.57); in the latter, $a = k$ in (2.60).  \hfill {\bf QED}

\medskip

Thus, the optimal investment strategy is same strategy as in the unconstrained case but truncated by $w$ if necessary.

Finally, we indicate how the results of this section are related to those under the same financial model (including the constraint of no borrowing) for an individual who maximizes her expected discounted utility of consumption, when the utility function exhibits constant relative risk aversion (CRRA).  We give the parameter values of the utility maximization problem that lead to {\it identical} consumption and investment strategies as we  found for the individual who minimizes her probability of lifetime ruin.  Young (2004) found a similar parallel in the unconstrained case. 

Specifically, the utility function is of the form

$$u(c) = {c^\eta \over \eta}, \quad \eta < 1, \eta \ne 0, \eqno(2.68)$$

\noindent and the corresponding value function is given by

$$V_0(w) = \sup_{\pi_0} {\bf E}_w \int_0^\infty e^{-\beta t} u(C_t) dt, \eqno(2.69)$$

\noindent in which $\beta$ is a personal discount factor, and the supremum is taken over admissible  constrained investment strategies as in the definition of $\psi_0$ in (2.3).

Let $\beta = \lambda + p$.  One can adapt the arguments in Fleming and Zariphopoulou (1991) to show that

\item{(i)}  If ${\mu - r \over \sigma^2} {1 \over a_r + 1} < 1$, then the investor who minimizes her probability of lifetime ruin behaves as a CRRA utility maximizer with $\eta = -a_r$, or equivalently, with CRRA = $1 + a_r$.

\item{(ii)}  If ${\mu - r \over \sigma^2} {1 \over a_r + 1} \ge 1$, then the investor who minimizes her probability of lifetime ruin behaves as a CRRA utility maximizer with $\eta = -k$, or equivalently, with CRRA = $1 + k$.

In other words, if we see an investor behaving as prescribed in this section, then we do not know whether she is minimizing her probability of lifetime ruin under proportional consumption or maximizing her utility of consumption under power utility.  It is interesting that we get this parallel when the personal discount factor $\beta = \lambda + p$.  The discount $\beta$ measures the individual's impatience.  The higher the hazard rate $\lambda$ or the proportion consumed $p$, the greater the impatience, which corresponds with our intuition.

\sect{3. Probability of Lifetime Ruin: Borrowing Rate Higher than Lending Rate}

In this section, we consider the problem of minimizing the probability of lifetime ruin when the individual is allowed to borrow money at a rate $b$ that is higher than the lending rate $r$.  In Section 3.1, we present the financial market.  In Section 3.2, we solve for the probability of ruin in the case for which consumption is constant.  Section 3.3 parallels Section 3.2 in the case for which consumption is proportional to wealth.  We show how this latter case is related to optimal investment and consumption in a Merton model with borrowing, as considered in Fleming and Zariphopoulou (1991).

\subsect{3.1. Financial Market}

We assume that the individual invests a non-negative amount in a riskless asset that earns interest at the constant rate $r > 0$.  If the individual borrows money, then she pays interest at the rate $b \ge r$.  Also, the individual invests in a risky asset whose price follows geometric Brownian motion, as given in (2.1), in which $\mu > b \ge r > 0$.

Let $W_t$ be the wealth at time $t$ of the individual, and let $\pi_{b,t}$ be the amount that the decision maker invests in the risky asset at that time.  We use a subscript $b$ to denote the fact that borrowing is allowed but only at a higher rate than the lending rate.  It follows that the amount invested in the riskless asset is $(W_t - \pi_{b,t})_+$, and the amount borrowed is $(\pi_{b,t} - W_t)_+$. Thus, wealth follows the process

$$\left\{ \eqalign{dW_t &= [r (W_t - \pi_{b, t})_+ - b (\pi_{b,t} - W_t)_+ + \mu \pi_{b,t}  - c(W_t)] dt + \sigma \pi_{b,t} dB_t, \cr
W_0 &= w.} \right. \eqno(3.1)$$

In this case, we denote the minimal probability of ruin by $\psi_b$ and define it similarly as in (2.3), except that an admissible investment strategy allows borrowing and short-selling.  For every $\alpha \in {\bf R}$, we associate a second-order differential operator ${\cal D}^{\alpha}$ with this minimization problem as follows:  For every open set $G \subset {\bf R}^+$ and for every $h \in C^2(G)$, we define the function ${\cal D}^{\alpha} h: G \rightarrow {\bf R}$ by

$${\cal D}^{\alpha} h(w) = \left[ bw + (\mu - b) \alpha - c(w) \right] h'(w) + {1 \over 2} \sigma^2 \alpha^2 h''(w) - \lambda h(w). \eqno(3.2)$$

\noindent Also, recall the definition of ${\cal L}^{\alpha}$ given in (2.4).  We use both differential operators in the following sections to characterize the minimum probability of ruin $\psi_b$ in the cases for which $c(w) = c$ and $c(w) = pw$, respectively.

\subsect{3.2. Constant Consumption}

In this section, we consider the case for which the consumption rate $c(w) = c$, a positive constant, as in Section 2.2.  As argued in that section, we have the probability of ruin equal to zero when wealth is at least $c/r$.  We begin with a verification theorem whose proof we omit because it is similar to the one of Theorem 2.1.

\th{3.1} {Suppose $h_b$ is a decreasing function from ${\bf R}^+$ to $[0, 1]$ and suppose $\alpha_b$ is a function from $[0, c/r]$ to $\bf R$ that satisfy the following conditions:
\item{$(i)$} $h_b \in C^2$ on $[0, c/r);$
\item{$(ii)$} ${\cal L}^\alpha h_b(w) \ge 0,$ for $w \in (0, c/r)$ and $\alpha \in {\bf R};$
\item{$(iii)$} If $0 \le \alpha_b(w) \le w,$ then ${\cal L}^{\alpha_b(w)} h_b(w) = 0;$
\item{$(iv)$} ${\cal D}^\alpha h_b(w) \ge 0,$ for $w \in (0, c/r)$ and $\alpha \in {\bf R};$
\item{$(v)$} If $\alpha_b(w) \ge w,$ then ${\cal D}^{\alpha_b(w)} h_b(w) = 0;$
\item{$(vi)$} $h_b(0) = 1,$ $h_b(w) = 0 = \alpha_b(w)$ for $w \ge c/r$.}

{\it Under the above conditions, the minimum probability of the lifetime ruin $\psi_b$ is given by}
$$\psi_b(w) = h_b(w), \quad w \ge 0, \eqno(3.3)$$

\noindent {\it and the optimal investment strategy in the risky asset $\pi^*_b$ is given by}
$$\pi^*_b(w) = \alpha_b(w), \quad w \in [0, c/r]. \eqno(3.4)$$

\medskip

As in Section 2.2, we demonstrate the steps of how to find $h_b$ and $\alpha_b$ that satisfy the hypotheses of Theorem 3.1.  Recall that in the case for which $b = r$, the amount invested in the riskless asset increases (linearly) from a negative amount when $w = 0$ (that is, the individual borrows money to invest in the risky asset) to a positive amount when $w = c/r$.

{\bf Ansatz:} We use this result to hypothesize that in the case with $\mu > b \ge r$, the amount invested in the riskless asset satisfies

$$w - \alpha_b(w) \cases{ < 0, & if $w < w_b$; \cr = 0, & if $w_b \le w \le w_l$; \cr > 0, & if $w_l < w \le c/r$;} \eqno(3.5)$$

\noindent for some constants $0 \le w_b \le w_l \le c/r$.  In addition, we hypothesize that $\alpha_b$ is continuous on $[0, c/r]$.  Finally, we suppose that $h_b$ is strictly decreasing on $[0, c/r]$.

\medskip

Under the ansatz, we, therefore, have three intervals of wealth to consider:  $(w_l, c/r]$, $[w_b, w_l]$, and $[0, w_b)$.  First, consider $(w_l, c/r]$, on which the individual invests a positive amount of money in the riskless asset at the rate $r$.  We have the following proposition concerning $h_b$ on this interval.  The proof is identical to that of Propositions 2.2, so we omit it.

\prop{3.2} {Under the hypothesis that $0 \le \alpha_b(w) < w$ on $(w_l, c/r],$ the function $h_b$ given by}

$$h_b(w) = \beta_b \left( 1 - {r \over c} w \right)^d, \eqno(3.6)$$

\noindent{\it with $d$ given in $(2.15)$ and $\beta_0 \ge \beta_b \ge 1$, satisfies the properties specified in Theorem 3.1.  Thus, on $(w_l, c/r]$, $h_b$ is a multiple of the probability of ruin $\psi$ in (2.14).  The corresponding $\alpha_b$ on this interval is given by}

$$\alpha_b(w) = {\mu - r \over \sigma^2} {1 \over d - 1} \left( {c \over r} - w \right),  \eqno(3.7)$$

\noindent {\it identical to $\pi^*$ in the case for which $b = r$ and to $\pi^*_0$ in the no-borrowing case.}

\medskip

\cor{3.3} {By the assumed continuity of $\alpha_b,$ it follows from $(3.7)$ that $w_l$ is given by $(2.33)$.}

\medskip

In other words, the lending level $w_l$ in the case with borrowing is identical to the lending level in the case for which $b = r$ and in the no-borrowing case.  Thus, we see the same myopia in this case that we saw in the no-borrowing case.

Next, consider the interval $[w_b, w_l]$, on which $\alpha_b(w) = w$.  We have the following proposition concerning $h_b$ on this interval that we state without proof.

\prop{3.4} {If $h_b$ exists as specified in Theorem 3.1, then under the hypothesis that $\alpha_b(w) = w$ on $[w_b, w_l]$, $h_b$ solves}

$$\lambda h_b = (\mu w - c) h_b' + {1 \over 2} \sigma^2 w^2 h_b'',  \eqno(3.8)$$

\noindent {\it with boundary condition}

$${h_b(w_l) \over h'_b(w_l)} = -{1 \over d} \left( {c \over r} - w_l \right). \eqno(3.9)$$

\medskip

As in Section 2.2, we similarly define $y$ on $[w_b, w_l]$ by

$$y(w) = {h_b(w) \over h'_b(w)}. \eqno(3.10)$$

\noindent Note that $y$ solves the differential equation (2.44) with boundary condition (3.9) and (3.10) at $w = w_l$.  Thus, the function $y$ in (3.10) is identical to the one in Section 2.2.  For that reason, we use the same letter $y$ to denote this function.  We use this function later in determining $h_b$ on $[0, w_b)$, including determining $w_b$ itself.

Finally, consider the interval $[0, w_b)$, on which we hypothesize that the individual borrows money at the rate $b$.  Before stating a proposition concerning $h_b$ on this interval, we outline a program whereby one can determine $h_b$ from its Legendre transform $\tilde h$.  By items (v) and (vi) of Theorem 3.1, the function $h_b$ solves

$$\lambda h_b = (bw - c) h'_b + \min_{\alpha}  \left[ (\mu - b) \alpha h'_b(w) + {1 \over 2} \sigma^2 \alpha^2 h''_b(w) \right], \eqno(3.11)$$

\noindent with boundary condition $h_b(0) = 1$.  Note that (3.11) is a fully-nonlinear differential equation, but we can linearize (3.11) by rewriting it in terms of the Legendre transform of $h_b$.

Define the Legendre transform $\tilde h$ of $h_b$ on $[0, w_b)$ by

$$\tilde h(v) = \min_w [h_b(w) + wv]. \eqno(3.12)$$

\noindent Note that we can recover $h_b$ from $\tilde h$ by

$$h_b(w) = \max_v [\tilde h(v) - wv]. \eqno(3.13)$$

\noindent The minimizing value of $w$ in (3.12) equals $I(-v) = \tilde h'(v)$, in which $I$ is the inverse function of $h'_b$.  Therefore, the maximizing value of $v$ in (3.13) equals $-h'_b(w)$.

Substitute $w = I(-v)$ in equation (3.11) to obtain

$$\lambda \tilde h(v) + (b - \lambda) v \tilde h'(v) - m_b v^2 \tilde h''(v) = cv, \eqno(3.14)$$

\noindent in which $m_b$ is defined as in (2.16) with $r$ replaced by $b$.  The general solution of (3.14) is

$$\tilde h(v) = D_1 v^{B_1} + D_2 v^{B_2} + {c \over b} v, \eqno(3.15)$$

\noindent in which $D_1$ and $D_2$ are constants to be determined, and $B_1$ and $B_2$ are the positive and negative roots, respectively of

$$-\lambda - (b - \lambda + m_b) B + m_b B^2 = 0.  \eqno(3.16)$$

\noindent Thus,

$$B_1 = {1 \over 2 m_b} \left[ (b - \lambda + m_b) + \sqrt{(b - \lambda + m_b)^2 + 4 \lambda m_b} \right] > 1, \eqno(3.17)$$

\noindent and

$$B_2 = {1 \over 2 m_b} \left[ (b - \lambda + m_b) - \sqrt{(b - \lambda + m_b)^2 + 4 \lambda m_b} \right] < 0. \eqno(3.18)$$

Define $v_0 = -h'_b(0)$ and $v_b = -h'_b(w_b)$; that is,

$$\tilde h'(v_0) = 0, \eqno(3.19)$$ and

$$\tilde h'(v_b) = w_b. \eqno(3.20)$$  Recall that we still do not know $w_b$.

From the definition of $\tilde h$ in (3.12), at $v = v_0$, we have

$$\tilde h(v_0) = 1. \eqno(3.21)$$

\noindent At $v = v_b$, we use the assumed continuity of $\alpha_b(w) = -{\mu - b \over \sigma^2} {h'_b(w) \over h''_b(w)}$ at $w = w_b$ to get

$$\tilde h'(v_b) = - {\mu - b \over \sigma^2} v_b \tilde h''(v_b). \eqno(3.22)$$

To summarize where we are at this point:  To determine $h_b$ on $[0, c/r]$, we need to know $\beta_b$, $w_b$, $v_b$, $v_0$, $D_1$, and $D_2$.  In what follows, we use expressions $(3.19) - (3.22)$, together with smoothness of $h_b$, to calculate these six unknowns.

For the moment, assume that we know $w_b$, and later we show how to obtain it from the function $y$.  Equations (3.20) and (3.22) imply, respectively,

$$D_1 B_1 v_b^{B_1} + D_2 B_2 v_b^{B_2} + {c \over b} v_b = w_b v_b, \eqno(3.23)$$ and

$$D_1 B_1 (B_1 - 1) v_b^{B_1} + D_2 B_2 (B_2 - 1) v_b^{B_2} = -w_b v_b {\sigma^2 \over \mu - b}. \eqno(3.24)$$

\noindent Solve (3.23) and (3.24) for $D_1$ and $D_2$ to get

$$D_1 = - {v_b^{1 - B_1} \over B_1 (B_1 - B_2)} \left[ w_b {\sigma^2 \over \mu - b} + \left( {c \over b} - w_b \right) (1 - B_2) \right], \eqno(3.25)$$ and

$$D_2 = - {v_b^{1 - B_2} \over B_2 (B_1 - B_2)} \left[ -w_b {\sigma^2 \over \mu - b} + \left( {c \over b} - w_b \right) (B_1 - 1) \right]. \eqno(3.26)$$

Substitute for $D_1$ and $D_2$ in equation (3.19), that is, in $-D_1 B_1 v_0^{B_1 - 1} - D_2 B_2 v_0^{B_2 - 1} = {c \over b}$, to get

$$\eqalign{& {1 \over B_1 - B_2} \left( {v_0 \over v_b} \right)^{B_1 - 1} \left[ w_b {\sigma^2 \over \mu - b} + \left( {c \over b} - w_b \right) (1 - B_2) \right] \cr
& + {1 \over B_1 - B_2} \left( {v_0 \over v_b} \right)^{B_2 - 1} \left[ -w_b {\sigma^2 \over \mu - b} + \left( {c \over b} - w_b \right) (B_1 - 1) \right] = {c \over b}.} \eqno(3.27)$$

\noindent We next show that (3.27) has a unique root $v_0/v_b > 1$.  First, note that if $v_0/v_b = 1$, then the left-hand-side of (3.27) equals $c/b - w_b < c/b$.  Next, note that as $v_0/v_b \rightarrow \infty$, the left-hand-side of (3.27) goes to infinity if

$$w_b {\sigma^2 \over \mu - b} + \left( {c \over b} - w_b \right) (1 - B_2) > 0, \eqno(3.28)$$

\noindent which we will show is true when we find $w_b$.  Finally, by taking the derivative of the left-hand-side of (3.27) with respect to $v_0/v_b > 1$, one can show that the left-hand-side of (3.27) is increasing if (3.28) holds.

Once we have the solution $v_0/v_b > 1$ of equation (3.27), we can substitute it into (3.21) via $D_1$ and $D_2$ and solve for $v_0$ as follows:

$$\eqalign{& - {1 \over B_1(B_1 - B_2)} \left( {v_0 \over v_b} \right)^{B_1 - 1} \left[ w_b {\sigma^2 \over \mu - b} + \left( {c \over b} - w_b \right) (1 - B_2) \right] \cr
& - {1 \over B_2(B_1 - B_2)} \left( {v_0 \over v_b} \right)^{B_2 - 1} \left[ -w_b {\sigma^2 \over \mu - b} + \left( {c \over b} - w_b \right) (B_1 - 1) \right] + {c \over b} = {1 \over v_0}.} \eqno(3.29)$$

\noindent Then, we can get $v_b$ from

$$v_b = {v_0 \over v_0/v_b}, \eqno(3.30)$$

\noindent and by substituting $v_b$ into (3.25) and (3.26), we have $D_1$ and $D_2$, respectively.

Now, it only remains for us to determine $w_b$ and verify inequality (3.28).  Once we have $w_b$, we calculate $h_b(w_b)$ by continuity of the solution of $h_b$ on $[0, w_b)$, from which we get $h_b$ on $[w_b, w_l]$ by numerically solving the boundary-value problem in (3.8) and (3.9).  In particular, we compute $h_b(w_l)$, from which $\beta_b$ follows.

\lem{3.5} {$w_b$ is the unique solution of}

$$y(w_b) = {\mu + b \over 2 \lambda} w_b - {c \over \lambda}. \eqno(3.31)$$

\medskip

\pf  By assumed smoothness of $h_b$, the differential equation (2.35) holds at $w = w_b$, with $h_b$ substituted for $h_0$; specifically,

$$\lambda h_b(w_b) = (\mu w_b - c) h_b'(w_b) + {1 \over 2} \sigma^2 w_b^2 h_b''(w_b).  \eqno(3.32)$$

\noindent Into equation (3.32), substitute 

$$w_b = -{\mu - b \over \sigma^2} {h'_b(w_b) \over h''_b(w_b)} \eqno(3.33)$$

\noindent to get

$$\lambda h_b(w_b) = (\mu w_b - c) h_b'(w_b) - {1 \over 2} (\mu - b) w_b h_b'(w_b). \eqno(3.34)$$

\noindent Recall that $y(w_b) = h_b(w_b)/h'_b(w_b)$; thus, (3.34) becomes (3.31).

Note that the right-hand-side of (3.31) is of the same form as $z$ given in (2.43) with $r$ replaced by $b$, so denote the function implied by the right-hand-side of (3.31) by $z_b$.  Thus, $w_b$ is defined similarly to $w_l$ because $y(w_l) = {\mu + r \over 2 \lambda} w_l - {c \over \lambda}$.

It remains for us to show that $y$ intersects $z_b$ uniquely.  First, suppose that $y$ intersects $z_b$ in two points $0 < w_b^* < w_b < w_l$ with $y'(w_b^*) \ge z'_b(w_b^*) = {\mu + b \over 2 \lambda}$ and $y'(w_b) \le z'_b(w_b) = {\mu + b \over 2 \lambda}$.  From the proof of Lemma 2.6, we have that

$$1 - {\mu + b \over 2 \lambda} \ge - {\mu - b \over \sigma^2 w_b^*} z_b(w_b^*). \eqno(3.35)$$

\noindent The right-hand-side of (3.35) is positive, so if ${\mu + b \over 2 \lambda} \ge 1$, then we have our contradiction.

Therefore, suppose ${\mu + b \over 2 \lambda} < 1$, from which it follows as in the proof of Lemma 2.6 that

$$w_b^* \ge {2c(\mu - b) \over \sigma^2(2 \lambda - \mu - b) + (\mu^2 - b^2)}.  \eqno(3.36)$$

\noindent On the other hand, $y'(w_b) \le z'_b(w_b)$ implies that

$$w_b \le {2c(\mu - b) \over \sigma^2(2 \lambda - \mu - b) + (\mu^2 - b^2)}, \eqno(3.37)$$

\noindent a contradiction to $w_b^* < w_b$.  Thus, if $y$ intersects $z_b$, it can do so in only one point.

Finally, we show that $y$ intersects $z_b$ for all $\mu > b > r$.  It is enough to show that $y'(0) = \mu/\lambda$.  From $y$'s differential equation (2.44), we obtain that $y(0) = -c/\lambda$.  Differentiate (2.44) and subsitute $w = 0$ to get $y'(0) = \mu/\lambda$.  Thus, we are done. \hfill {\bf QED}

\medskip

In a corollary to Lemma 3.5, we assert that inequality (3.28) holds.

\cor{3.6} {Inequality $(3.28)$ holds.}

\medskip

\pf It is enough to show that $w_b < c/b$.  $z_b$ intersects the horizontal axis at $w = 2c/(\mu + b) > w_b$.  Now, $2c/(\mu + b) < c/b$ because $\mu > b$; thus, $w_b < c/b$.  \hfill {\bf QED}

\medskip

We summarize the solution of $h_b$ on $(0, w_b)$ in the following proposition.

\prop{3.7} {On $(0, w_b)$, in which $w_b$ is given in Lemma 3.5, first solve $(3.27)$ for $v_0/v_b > 1$, then solve $(3.29)$ and $(3.30)$ for $v_0$ and $v_b$, respectively.  Then, use $v_b$ and $w_b$ in the expressions for $D_1$ and $D_2$, namely $(3.25)$ and $(3.26),$ respectively.  Next, determine $\tilde h$ in $(3.15)$. Finally, obtain $h_b$ from $\tilde h$ via the transform in $(3.13)$.}

\medskip

Next, we show that if $h_b$ solves (3.8) and the corresponding boundary conditions, then the optimal investment in the risky asset is equal to the current wealth.  We state this formally in the following proposition.

\prop{3.8} {Suppose $h_b$ solves $(3.8), (3.9),$ and $h_b(w_b+) = h_b(w_b-)$, in which $w_b$ is given in Lemma 3.5; then,}

$$\eqalign{\arg \min & \left( \min_{\alpha \ge w} \left[ (\mu - b) \alpha h'_b(w) + {1 \over 2} \sigma^2 \alpha^2 h''_b(w) \right], \right. \cr
& \left. \quad \min_{\alpha \le w} \left[ (\mu - r) \alpha h'_b(w) + {1 \over 2} \sigma^2 \alpha^2 h''_b(w) \right] \right) = w.} \eqno(3.38)$$

\medskip

\pf Consider the parabolas $f_1$ and $f_2$ defined by

$$f_1(\pi) = (\mu - b) \pi h'_b(w) + {1 \over 2} \sigma^2 \pi^2 h''_b(w), \eqno(3.39)$$
 
\noindent and

$$f_2(\pi) = (\mu - r) \pi h'_b(w) + {1 \over 2} \sigma^2 \pi^2 h''_b(w). \eqno(3.40)$$

\noindent  The minimum of $f_1$ for $\pi \ge w$ occurs at $\pi = w$ if $f'_1(w) \ge 0$, which is true if $y < z_b$ on $(w_b, w_l)$.  From the proof of Lemma 3.5, this latter condition is true.  The minimum of $f_2$ for $\pi \le w$ occurs at $\pi = w$ if $f'_2(w) \le 0$, which is true if $y > z$ on $(w_b, w_l)$.  From Lemma 2.6, this latter condition is true.  Therefore, the minimum of the expression on the left-hand-side of (3.38) occurs at $\pi = w$.  \hfill {\bf QED}

\medskip

Parallel to Proposition 3.8, we have the following proposition that asserts that the optimal investment in the risky asset on $(0, w_b)$ exceeds the current wealth; therefore, the individual does borrow money when wealth is less than $w_b$.

\prop{3.9} {On $(0, w_b)$, the optimal policy $\alpha_b$ satisfies}

$$\alpha_b(w) > w. \eqno(3.41)$$

\medskip

\pf $\alpha_b(w) > w$ if and only if

$$- {\mu - b \over \sigma^2} v \tilde h''(v) > \tilde h'(v). \eqno(3.42)$$

\noindent By substituting for $\tilde h$ in (3.42) and by using the expressions in (3.25) and (3.26) for $D_1$ and $D_2$, respectively, we obtain the following equivalent inequality.

$$\eqalign{& {1 \over B_1 - B_2} \left( {v \over v_b} \right)^{B_1 - 1} \left( {\mu - b \over \sigma^2} (B_1 - 1) + 1 \right) \left( w_b {\sigma^2 \over \mu - b} + \left( {c \over b} - w_b \right) (1 - B_2) \right) \cr
& + {1 \over B_1 - B_2} \left( {v \over v_b} \right)^{B_2 - 1} \left( {\mu - b \over \sigma^2} (B_2 - 1) + 1 \right) \left( -w_b {\sigma^2 \over \mu - b} + \left( {c \over b} - w_b \right) (B_1 - 1) \right) > {c \over b}.} \eqno(3.43)$$

At $v = v_b$, the left-hand-side of (3.43) equals $c/b$.  Therefore, if the left-hand-side of (3.43) increases with respect to $v/v_b$, then we are done.  The derivative of the left-hand-side of (3.43) is positive if and only if inequality (3.28) holds (which we know to be true from Corollary 3.6) and if

$$\eqalign{& {B_1 - 1 \over B_1 - B_2} \left( {\mu - b \over \sigma^2} (B_1 - 1) + 1 \right) \left( w_b {\sigma^2 \over \mu - b} + \left( {c \over b} - w_b \right) (1 - B_2) \right) \cr
& + {B_2 - 1 \over B_1 - B_2} \left( {\mu - b \over \sigma^2} (B_2 - 1) + 1 \right) \left( -w_b {\sigma^2 \over \mu - b} + \left( {c \over b} - w_b \right) (B_1 - 1) \right) > 0.} \eqno(3.44)$$

\noindent By simplifying inequality (3.44), one can show that it is equivalent to

$$w_b (\sigma^2 (\mu + b - 2 \lambda) - (\mu - b)^2) + 2 b (\mu - b) \left( {c \over b} - w_b \right) > 0. \eqno(3.45)$$

In order to show inequality (3.45), we consider two cases:  $(\mu + b)/(2 \lambda) \ge 1$, and $(\mu + b)/(2 \lambda) < 1$.  First, suppose $(\mu + b)/(2 \lambda) \ge 1$.  If $\sigma^2 (\mu + b - 2 \lambda) - (\mu - b)^2 \ge 0$, then we are done because $c/b - w_b > 0$.  On the other hand, if $\sigma^2 (\mu + b - 2 \lambda) - (\mu - b)^2 < 0$, then inequality (3.45) holds if

$${2c \over \mu + b} (\sigma^2 (\mu + b - 2 \lambda) - (\mu - b)^2) + 2 b (\mu - b) \left( {c \over b} - w_b \right) \ge 0 \eqno(3.46)$$

\noindent because we replaced $w_b$ with something larger.  Now, inequality (3.46) reduces to

$$\sigma^2 (\mu + b - 2 \lambda) \ge 0, \eqno(3.47)$$

\noindent which is true by assumption.

Next, suppose $(\mu + b)/(2 \lambda) < 1$.  Then, $y(w_b) = z_b(w_b)$ and $y'(w_b) > z'_b(w_b)$ imply that

$$w_b < {2c(\mu - b) \over \sigma^2(2 \lambda - \mu - b) + (\mu^2 - b^2)}, \eqno(3.48)$$

\noindent as in inequality (3.37).  Thus, inequality (3.45) holds if

$$\eqalign{&{2c(\mu - b) \over \sigma^2(2 \lambda - \mu - b) + (\mu^2 - b^2)} (\sigma^2 (\mu + b - 2 \lambda) - (\mu - b)^2) \cr
& \quad + 2 b (\mu - b) \left( {c \over b} - {2c(\mu - b) \over \sigma^2(2 \lambda - \mu - b) + (\mu^2 - b^2)} \right) \ge 0} \eqno(3.49)$$

\noindent because we replaced $w_b$ with something larger.  The left-hand-side of inequality (3.49) reduces to 0, so we are done.  \hfill {\bf QED}

\medskip

We have the following theorem that summarizes the results of this section.

\th{3.10} {For constant consumption $c(w) = c$, the borrowing level $w_b$ is given by Lemma 3.5, and the lending level is (2.33).  The minimum probability of ruin $\psi_b$ is given by}

$$\psi_b(w) = h_b(w), \quad 0 \le w < w_b, \eqno(3.50)$$

\noindent {\it in which $h_b$ is obtained as described in Proposition 3.7};

$$\psi_b(w) = h_b(w), \quad w_b \le w \le w_l, \eqno(3.51)$$

\noindent {\it in which $h_b$ is obtained from Proposition 3.4 and from continuity of $h_b$ at $w_b;$ and}

$$\psi_b(w) = \beta_b \left( 1 - {r \over c} w \right)^d, \quad w_l < w \le c/r, \eqno(3.52)$$

\noindent {\it in which $\beta_b = h_b(w_l)(1 - rw_l/c)^{-d}$.  Finally, $\psi_b(w) = 0$ for $w > c/r$.}

{\it The optimal investment strategy is given by}

$$\pi^*_b(w) = \cases{-{\mu - b \over \sigma^2} v \tilde h''(v), & if $0 \le w < w_b$, for $\tilde h'(v) = w;$ \cr w, & if $w_b \le w \le w_l;$ \cr {\mu - r \over \sigma^2} {1 \over d - 1} \left( {c \over r} - w \right), & if $w_l < w \le c/r$,} \eqno(3.53)$$

\noindent{\it in which $\tilde h$ is given by $(3.12)$ on $[0, w_b)$.}

\medskip

During the remainder of this section, we examine the limit of the optimal investment strategy as the borrowing rate $b$ approaches the drift on the risky asset $\mu$ from the left.  We show that the optimal investment strategy approaches that in the no-borrowing case if $\mu \le \lambda$; however, if $\mu > \lambda$, then the amount invested in the risky asset is arbitrarily large for wealth near zero.  Thus, instead of reducing leveraging by imposing a higher borrowing rate, leveraging actually {\it increases} if $\mu > \lambda$.

We begin with the following lemma.

\lem{3.11} {The probability of lifetime ruin $\psi_b$ is convex on $(0, c/r)$.  In addition,
\item{$(a)$} If $\mu \le \lambda,$ then $\psi_0$ is convex on $(0, c/r)$.
\item{$(b)$} If $\mu > \lambda,$ then $\psi_0$ is convex on $(w_\mu, c/r)$ but is concave on $(0, w_\mu)$, in which \hfill \break $\lim_{b \rightarrow \mu-} w_b = w_\mu > 0$, with $w_b$ given in Lemma 3.5.}

\medskip

\pf We prove this lemma by considering $\lim_{b \rightarrow \mu-} w_b = w_\mu$.  We show that $w_\mu = 0$ when $\psi_0$ has no inflection point, that is, when $\psi_0$ is convex on $(0, c/r)$.  This occurs when $\mu \le \lambda$.  Similarly, we show that $w_\mu > 0$ when $\psi_0$ changes concavity; $w_\mu$ is the inflection point of $\psi_0$ in this case.  This occurs when $\mu > \lambda$.

Recall that for $\mu > b \ge r$, $y$ intersects $z_b$ at $w = 0$ and at $w = w_b \in (0, w_l]$ with $y > z_b$ on $(0, w_b)$.  By continuity of $y$, if we take the limit as $b$ approaches $\mu$ from the left, then $y$ will intersect $z_\mu$ at $w = 0$ and at $w = w_\mu \ge 0$ with $y > z_\mu$ on $(0, w_\mu)$, in which $z_\mu(w) = (\mu/\lambda)w - c/\lambda$.  Note that $(0, w_\mu)$ is empty if $w_\mu = 0$.

Calculate $y''(0)$ by differentiating equation (2.44) twice and evaluating the result at $w = 0$:  $y''(0) = -\lambda(\mu - \lambda)/c^2$.  Recall that $y'(0) = \mu/\lambda = z'_\mu(0)$.  If $\mu \le \lambda$, then $y''(0) \le 0$ from which it follows that $y \le z_\mu$ on $(0, \delta)$ for some $\delta > 0$.  Thus, $w_\mu = 0$ in this case because $y$ does not intersection $z_\mu$.

If $\mu > \lambda$, then $y''(0) > 0$ and $\psi''_0(0) < 0$.  From equation (2.35), note that $\psi''_0$ is proportional to $z_\mu - y$ on $(0, w_l]$.  Thus, $y$ intersects $z_\mu$ at $w = 0$ and $w = w_\mu > 0$ with $y > z_\mu$ on $(0, w_\mu)$ and $y < z_\mu$ on $(w_\mu, w_l)$.  Equivalently, $\psi''_0 < 0$ on $(0, w_\mu)$, and $\psi''_0 > 0$ on $(w_\mu, w_l)$, from which it follows that $w_\mu$ is the inflection point of $\psi_0$.  \hfill {\bf QED}

\th{3.12} {Let $w_b$ be as given in Lemma 3.5.  As $b \rightarrow r+$, $\psi_b$ and $\pi^*_b$ converge to the probability of ruin and investment strategy, $\psi$ in (2.14) and $\pi^*$ in (2.17), respectively.
\item{$(a)$} If $\mu \le \lambda,$ then as $b \rightarrow \mu-$, $\psi_b$ and $\pi^*_b$ converge to $\psi_0$ and $\pi^*_0$, respectively.
\item{$(b)$} If $\mu > \lambda,$ then as $b \rightarrow \mu-$, $\psi_b$ does not converge to $\psi_0$.  Moreover, $\pi^*_b(w)$ becomes arbitrarily large for wealth near 0 as $b \rightarrow \mu-$.}  

\medskip

\pf Because $\lim_{b \rightarrow r+} w_b = w_l$, $\psi_b$ and $\pi^*_b$ converge to the probability of ruin and investment strategy, $\psi$ and $\pi^*$, respectively, as $b \rightarrow r+$.

If $\mu \le \lambda$, then the proof of Lemma 3.11 gives us that  $\lim_{b \rightarrow \mu-} w_b = 0$.  It follows easily that $\psi_b$ and $\pi^*_b$ converge to $\psi_0$ and $\pi^*_0$, respectively, in this case.  Next, from financial arguments, we have that $\psi_b$ increases as $b \rightarrow \mu-$; thus, $\lim_{b \rightarrow \mu-} \psi_b(w) = \sup_{r \le b < \mu} \psi_b(w)$. This limit is convex, since the supremum of convex functions is convex (Rockafellar, 1970, Theorem 5.5).  If $\mu > \lambda$, then $\psi_0$ is not convex; therefore, $\psi_b$ does not converge to $\psi_0$.

It remains for us to show that if $\mu > \lambda$, then $\pi^*_b(w)$ becomes arbitrarily large for wealth $w$ near $0$.  To this end, consider the expression for the optimal investment strategy in terms of the Legendre transform $\tilde h$ in (3.15): $\pi^*_b(w) = -{\mu - b \over \sigma^2} v \tilde h''(v)$, in which $v = -h'_b(w)$.  In terms of the dual variable $v$, the optimal investment strategy on $[v_b, v_0]$, after substituting for $D_1$ and $D_2$ from (3.25) and (3.26), respectively, is given by

$$\eqalign{& {\mu - b \over \sigma^2} {B_1 - 1 \over B_1 - B_2} \left( {v \over v_b} \right)^{B_1 - 1} \left( w_b {\sigma^2 \over \mu - b} + \left( {c \over b} - w_b \right) (1 - B_2) \right) \cr
& + {\mu - b \over \sigma^2} {B_2 - 1 \over B_1 - B_2} \left( {v \over v_b} \right)^{B_2 - 1} \left( -w_b {\sigma^2 \over \mu - b} + \left( {c \over b} - w_b \right) (B_1 - 1) \right).} \eqno(3.54)$$

\noindent Let $v = v_0$ in this expression, which corresponds to wealth equal to 0.  From equation (3.27), this expression becomes

$$\eqalign{& {\mu - b \over \sigma^2} {B_1 - 1 \over B_1 - B_2} \left( {v_0 \over v_b} \right)^{B_1 - 1} \left( w_b {\sigma^2 \over \mu - b} + \left( {c \over b} - w_b \right) (1 - B_2) \right) \cr
& + {\mu - b \over \sigma^2} (B_2 - 1) \left[ {c \over b} - {1 \over B_1 - B_2} \left( {v_0 \over v_b} \right)^{B_1 - 1} \left( w_b {\sigma^2 \over \mu - b} + \left( {c \over b} - w_b \right) (1 - B_2) \right) \right] \cr
& = {\mu - b \over \sigma^2} (B_2 - 1) {c \over b} + \left( {v_0 \over v_b} \right)^{B_1 - 1} \left( w_b + {\mu - b \over \sigma^2} \left( {c \over b} - w_b \right) (1 - B_2) \right).} \eqno(3.55)$$

\noindent Now, as $b$ approaches $\mu$ from the left, ${\mu - b \over \sigma^2} B_2$ goes to 0; however, $v_0/v_b$ is greater than 1 and $B_1$ goes to infinity.  Thus, the expression in (3.55) goes to infinity, and we are done.  \hfill {\bf QED}

\medskip

It is quite interesting that the limiting behavior of $\psi_b$ and $\pi^*_b$ depends on the relative values of the drift on the risky asset and the hazard rate.  If the drift on the risky asset is less than the hazard rate, then as the borrowing rate approaches the drift on the risky asset, the individual borrows less and less.  The probability of dying is great enough that she does not need to take on the risk of borrowing in order to avoid running out of money before dying.  On the other hand, if the drift on the risky asset is greater than the hazard rate, then because she has to pay more interest, she gets less value from her leverage than at lower interest rates and therefore borrows more and more.

In Section 4, we present numerical examples that demonstrate the results of this section.

\subsect{3.3. Proportional Consumption}

In this section, we consider the case for which the consumption rate is proportional to wealth $c(w) = pw$, in which $r \le b < \min(\mu, p)$.  Note that if we specify that ruin occurs if wealth reaches 0, then the individual with this consumption function does not ruin.  Therefore, for this case, as in Section 2.3, we say that ruin occurs when wealth reaches $w_0 > 0$; that is, let $\tau_0$ in definition (2.3) be the (first) time that wealth reaches $w_0$.

We have the following verification theorem whose proof we omit because it is similar to the proof for Theorem 2.1.

\th{3.13} {Suppose $g_b$ is a strictly decreasing function from ${\bf R}^+$ to $[0, 1]$ and suppose $\gamma_b$ is a function from ${\bf R}^+$ to ${\bf R}$ that satisfy the following conditions:
\item{$(i)$} $g_b \in C^2 \left( {\bf R}^+ \right);$
\item{$(ii)$} ${\cal L}^\gamma g_b(w) \ge 0,$ for $w \ge 0$ and $\gamma \in {\bf R};$
\item{$(iii)$} If $\gamma_b(w) \le w,$ then ${\cal L}^{\gamma_b(w)} g_b(w) = 0;$
\item{$(iv)$} ${\cal D}^\gamma g_b(w) \ge 0,$ for $w \ge 0$ and $\gamma \in {\bf R};$
\item{$(v)$} If $\gamma_b(w) \ge w,$ then ${\cal D}^{\gamma_b(w)} g_b(w) = 0;$
\item{$(vi)$} $g_b(0) = 1,$ and $\lim_{w \rightarrow \infty} g_b(w) = 0$.}

{\it Under the above conditions, the minimum probability of the lifetime ruin $\psi_0$ is given by}
$$\psi_b(w) = g_b(w), \quad w \ge 0, \eqno(3.56)$$

\noindent {\it and the optimal investment strategy in the risky asset $\pi^*_0$ is given by}
$$\pi^*_b(w) = \gamma_b(w), \quad w \ge 0. \eqno(3.57)$$

\medskip

As in the proof of Theorem 2.9, we can show that $g_b$ is given by

$$g_b(w) = \left( {w \over w_0} \right)^{-a},  \eqno(3.58)$$

\noindent for some $a > 0$.  From items (ii) through (v) of Theorem 3.13, we have that $a$ solves

$$\eqalign{\lambda &= \min \left\{ a (p - r) + a \min_{0 \le \gamma \le w} \left[ - (\mu - r) {\gamma \over w} + {1 \over 2} \sigma^2 (a + 1) {\gamma^2 \over w^2}  \right], \right. \cr
& \qquad \qquad \left. a (p - b) + a \min_{\gamma \ge w} \left[ - (\mu - b) {\gamma \over w} + {1 \over 2} \sigma^2 (a + 1) {\gamma^2 \over w^2}  \right] \right\} .} \eqno(3.59)$$

The function $f_1$ given by $f_1(\gamma) = -(\mu - r) {\gamma \over w} + {1 \over 2} \sigma^2 (a + 1) {\gamma^2 \over w^2}$ is minimized on $0 \le \gamma \le w$ by $\gamma = w$ if ${\mu - r \over \sigma^2} {1 \over a + 1} \ge 1$; otherwise, $f_1$ is minimized by $\gamma = {\mu - r \over \sigma^2} {w \over a + 1}$.

Similarly, the function $f_2$ given by $f_2(\gamma) = -(\mu - b) {\gamma \over w} + {1 \over 2} \sigma^2 (a + 1) {\gamma^2 \over w^2}$ is minimized on $\gamma \ge w$ by $\gamma = w$ if ${\mu - b \over \sigma^2} {1 \over a + 1} < 1$; otherwise, $f_2$ is minimized by $\gamma = {\mu - b \over \sigma^2} {w \over a + 1}$.

We have three cases:  (1) ${\mu - r \over \sigma^2} {1 \over a + 1} < 1$; (2) ${\mu - b \over \sigma^2} {1 \over a + 1} \ge 1$; and (3) ${\mu - r \over \sigma^2} {1 \over a + 1} \ge 1$ and ${\mu - b \over \sigma^2} {1 \over a + 1} < 1$.  We consider each in turn.

\medskip

\noindent {\bf Case 1:} ${\mu - r \over \sigma^2} {1 \over a + 1} < 1$.

Because $b > r$, we also have ${\mu - b \over \sigma^2} {1 \over a + 1} < 1$; thus, equation (3.59) becomes

$$\eqalign{\lambda &= \min \left[ a(p - r) + a f_1 \left( {\mu - r \over \sigma^2} {w \over a + 1} \right), d(p - b) + d f_2(w) \right] \cr
&= a(p - r) - {a m \over a + 1},} \eqno(3.60)$$

\noindent from which it follows that $a = a_r$ in (2.57).

\medskip

\noindent {\bf Case 2:} ${\mu - b \over \sigma^2} {1 \over a + 1} \ge 1$.

Because $b > r$, we also have ${\mu - r \over \sigma^2} {1 \over a + 1} \ge 1$; thus, equation (3.59) becomes

$$\eqalign{\lambda &= \min \left[ a(p - r) + a f_1(w), d(p - b) + d f_2 \left( {\mu - b \over \sigma^2} {w \over a + 1} \right) \right] \cr
&= a(p - b) - {a m_b \over a + 1},} \eqno(3.61)$$

\noindent from which it follows that $a = a_b$ in which $a_b$ is given by (2.57) with $r$ and $m$ replaced by $b$ and $m_b$, respectively.

\medskip

\noindent {\bf Case 3:} ${\mu - r \over \sigma^2} {1 \over a + 1} \ge 1$ and ${\mu - b \over \sigma^2} {1 \over a + 1} < 1$.

In this case, equation (3.59) becomes

$$\eqalign{\lambda &= \min \left[ a(p - r) + a f_1(w), d(p - b) + d f_2(w) \right] \cr
&= a \left(p - \mu + {1 \over 2} \sigma^2 (a + 1) \right),} \eqno(3.62)$$

\noindent from which it follows that $a = k$ in (2.60).

\medskip

We summarize the results of this section in the following theorem.

\th{3.14} {The minimum probability of ruin is given by}

$$\psi_b(w) =  \left( {w \over w_0} \right)^{-a}, \eqno(3.63)$$

\noindent {\it in which $a > 0$ is given by}

$$a = \cases{a_r, & if ${\mu - r \over \sigma^2} {1 \over a_r + 1} < 1;$ \cr
k, & if ${\mu - b \over \sigma^2} {1 \over k + 1} < 1 \le {\mu - r \over \sigma^2} {1 \over k + 1};$ \cr
a_b, & if ${\mu - b \over \sigma^2} {1 \over a_b + 1} \ge 1.$} \eqno(3.64)$$

\noindent {\it The corresponding constrained optimal investment strategy is given by}

$$\pi^*_b(w) = \cases{{\mu - r \over \sigma^2} {w \over a_r + 1}, & if ${\mu - r \over \sigma^2} {1 \over a_r + 1} < 1;$ \cr
w, & if ${\mu - b \over \sigma^2} {1 \over k + 1} < 1 \le {\mu - r \over \sigma^2} {1 \over k + 1};$ \cr
{\mu - b \over \sigma^2} {w \over a_b + 1}, & if ${\mu - b \over \sigma^2} {1 \over a_b + 1} \ge 1.$}  \eqno(3.65)$$

\medskip

We have the following corollary of Theorem 3.14 that gives us the limit of $\psi_b$ and $\pi^*_b$ as $b$ approaches $\min(\mu, p)$.

\cor{3.15} {From $\lim_{b \rightarrow r+} a_b = a_r$, it follows that $\psi_b$ and $\pi^*_b$ converge to $\psi$ in (2.56) and $\pi^*$ and (2.58), respectively.  From $\lim_{b \rightarrow \min(\mu, p)-} {\mu - b \over \sigma^2} {1 \over a_b + 1} = 0$, it follows that}

$$\lim_{b \rightarrow \min(\mu, p)-} \psi_b(w) = \psi_0(w), \quad w \ge 0, \eqno(3.66)$$

\noindent {\it and}

$$\lim_{b \rightarrow \min(\mu, p)-} \pi^*_b(w) = \pi^*_0(w), \quad w \ge 0. \eqno(3.67)$$

\noindent {\it In words, the $b = r$ case is the limit of the borrowing case as the borrowing rate $b$ approaches the rate of return $r$ on the riskless asset.  Also, the no-borrowing case is the limit of the borrowing case as the borrowing rate approaches the minimum of the proportion consumed and the drift on the risky asset.}

\medskip

Finally, we indicate how the results of this section are related to those under the same financial model for an individual who maximizes her expected discounted utility of consumption, when the utility function exhibits CRRA.  We give the parameter values of the utility maximization problem that lead to {\it identical} consumption and investment strategies as we  found for the individual who minimizes her probability of lifetime ruin.

As in Section 2.3, let $\beta = \lambda + p$.  From Fleming and Zariphopoulou (1991), we can show that

\item{(i)}  If ${\mu - r \over \sigma^2} {1 \over a_r + 1} < 1$, then the investor who minimizes her probability of lifetime ruin behaves as a CRRA utility maximizer with $\eta = -a_r$, or equivalently, with CRRA = $1 + a_r$.

\item{(ii)}  If ${\mu - b \over \sigma^2} {1 \over k + 1} < 1 \le {\mu - r \over \sigma^2} {1 \over k + 1}$, then the investor who minimizes her probability of lifetime ruin behaves as a CRRA utility maximizer with $\eta = -k$, or equivalently, with CRRA = $1 + k$.

\item{(iii)}  If ${\mu - b \over \sigma^2} {1 \over a_b + 1} \ge 1$, then the investor who minimizes her probability of lifetime ruin behaves as a CRRA utility maximizer with $\eta = -a_b$, or equivalently, with CRRA = $1 + a_b$.

\sect{4. Numerical Example}

In this section, we present a numerical example to demonstrate the results in Sections 2.2 and 3.2.  Because the probabilities of ruin and optimal investment strategies are explicitly or implicitly analytical, they are easy to implement with mathematical software.

Assume that real consumption is given by $c(w) = c$; that is, consumption is a constant (real) dollar rate. By giving the riskless rate $r$, the borrowing rate $b$, and the drift of the risky asset $\mu$ in real terms, we thereby model constant consumption in terms of real dollars. Assume the following parameter values:

\item{$\bullet$} $\lambda = 0.04$; the hazard rate is constant such that the expected future lifetime is 25 years.

\item{$\bullet$} $r = 0.02$; the riskless rate of return is 2\% over inflation.

\item{$\bullet$} $b = 0.04$; the borrowing rate is 4\% over inflation.

\item{$\bullet$} $\mu = 0.06$; the risky asset's drift is 6\% over inflation.

\item{$\bullet$} $\sigma = 0.20$; the risky asset's volatility is 20\%.

\item{$\bullet$} $c = 1$; the individual consumes one unit of wealth per year.

\medskip

In this example, the lending level $w_l = 14.64$ and the borrowing level $w_b = 10.62$.  Thus, in the no-borrowing case, the individual will invest $\pi^*_0(w) = w$ in the risky asset when $w \le 14.64$ and will invest less than her wealth when $w > 14.64$.  In the case for which the individual is allowed to borrow at the real rate of 4\%, the individual will borrow money when $w < 10.62$, will invest exactly her wealth in the risky asset when $10.62 \le w \le 14.64$, and will invest less than her wealth when $w > 14.64$.  See Figure 1 for a graph of the functions $y$, $z$, and $z_b$.  Note that $y > z$ for $0 < w < 14.64$, $y > z_b$ for $0 < w < 10.62$, and $y < z_b$ for $10.62 < w < 14.64$, as expected from Lemma 3.5.

\medskip

\centerline{\bf Figure 1 about here}
   
\medskip

See Figure 2 for a graph of the functions $\psi$, $\psi_0$, and $\psi_b$, in which $\psi$ is the probability of lifetime ruin in the unconstrained model for which $b = r$.  Note that $\psi \le \psi_b \le \psi_0$ on $[0, c/r]$, as expected.  $\lim_{b \rightarrow \mu-} w_b = w_\mu = 7.39$ in this case; therefore, $\psi_0$ is not convex for $w < 7.39$. 

\medskip

\centerline{\bf Figure 2 about here}

\medskip

See Figure 3 for a graph of the functions $\zeta$, $\zeta_0$, and $\zeta_b$, the amount invested in the riskless asset under the model for which $b = r$, under the no-borrowing model, and under the model for which $b = 0.04 > r$, respectively.  Specifically, $\zeta(w) =  w - \pi^*(w)$, and similarly for $\zeta_0(w)$ and $\zeta_b(w)$.  Recall that $\zeta(w) < 0$ means that the individual borrows $\vert \zeta(w) \vert$ at rate $r$, while $\zeta_b(w) < 0$ means that the individual borrows $\vert \zeta_b(w) \vert$ at rate $b$.

For $w > 14.64$, the three graphs are identical.  For $w \le 14.64$, $\zeta_0 \equiv 0$, so it is difficult to see; similarly, for $\zeta_b$ on the interval $[10.62, 14.64]$.  Notice that for $w < 10.62$, the graph of the amount borrowed under the $b = 0.04$ case is steeper than under the case for which $b = r = 0.02$.  For wealth near zero, the individual borrows more money when the borrowing rate is higher because she does not derive as much benefit from the leveraging when she pays more interest on the amount she borrows.

\medskip

\centerline{\bf Figure 3 about here}

\medskip

In Figure 4, we plot the graphs of $\zeta_b$ for three values of $b$: 0.04, 0.055, and 0.059.  Note that $\pi^*_b$ diverges for wealth near zero, as we expect from Theorem 3.12.

\medskip

\centerline{\bf Figure 4 about here}

\sect{5. Summary and Future Research}

In this paper, we consider an investment decision problem for an individual who seeks to minimize the probability of  outliving her wealth under borrowing constraints.  We consider two market models with different borrowing constraints. In the first model, the individual is prohibited from borrowing, whereas in the second model, the individual borrows at a rate that is higher than the rate earned on a (positive) investment in the riskless asset. The individual either consumes at a constant (real) dollar rate, or she consumes a constant proportion of her wealth. In each case, we find the minimum probability of lifetime ruin together with the optimal investment strategy.

Under the no-borrowing constraint, when the consumption function is constant, the optimal investment  strategy turns out to be a truncated version of the optimal investment strategy in the unconstrained case. The wealth level above which the individual invests a positive amount in the riskless asset, ``the lending level," turns out to be the same as it is under the unconstrained case. Below that level, the optimal investment strategy prescribes to put all the wealth in the risky asset. 

In the market model with different borrowing and lending rates, when the consumption function is constant, the individual will borrow money to invest in the risky asset when the current wealth is lower than some ``borrowing  level," will invest money in the riskless asset when her current wealth is higher than some ``lending level" (which is greater than the borrowing level and is equal to the lending level in the case for which $b = r$), and will invest all her money in the risky asset if her current wealth is between the borrowing level and the lending level.

Somewhat surprisingly, the model with a higher borrowing rate does not always converge to the model with the no-borrowing constraint as the borrowing rate goes to the drift of the risky asset. It turns out that the convergence holds only if the hazard rate of the individual is greater than or equal to the drift of the risky asset.  On the other hand, when the hazard rate is lower than the drift of the risky asset, the risk of running out of money before dying is great enough for low wealth that the individual borrows increasingly more money as the borrowing rate increases. 

In both of the market models, when the  consumption rate is proportional to the wealth, the individual who minimizes the probability of lifetime ruin behaves like an individual who maximizes her discounted utility of consumption under the same borrowing constraint when the utility function exhibits constant relative risk aversion (CRRA).

When we include life annuities in the financial market model, the probability of ruin decreases since the price of an annuity that guarantees the fixed consumption rate is less than the wealth over which the lifetime ruin probability is zero (``the safe level") (Milevsky, Moore, and Young, 2004). Therefore, in future work, we will extend our market model to include life annuities and determine the optimal investment allocations among risky and riskless assets and life annuities for an individual who minimizes her probability of lifetime ruin under borrowing constraints assumed here.  We also would like to extend our results to the more realistic case in which the hazard rate is time varying and apply the optimal stopping formulation developed in Milevsky, Moore, and Young (2004) to solve the problem given in this paper with borrowing constraints.

Also, comonotonicity and its applications to finance and insurance have been extensively studied by Goovaerts, Dhaene, Kaas, Denuit, and their co-workers; for example, see Dhaene et al. (2005) and Vanduffel, Dhaene, and Goovaerts (2005).  Amongst other applications, they propose comonotonic approximations to tackle multi-period optimal portfolio selection problems under very flexible deterministic saving and consumption patterns.  Their techniques might be useful to apply in our setting under more realistic modeling assumptions, such as random interest rates, more general risky asset price processes, and random consumption processes.  

\medskip

\centerline{\bf Acknowledgements} \medskip

We thank Moshe Milevsky, Kristen S. Moore, and S. David Promislow for their helpful comments.

\sect{References}

\smallskip \noindent \hangindent = 20 pt Browne, S. (1995), Optimal investment policies for a firm with a random risk process: Exponential utility and minimizing the probability of ruin, {\it Mathematics of Operations Research}, 20 (4): 937-958.

\smallskip \noindent \hangindent = 20 pt Browne, S. (1997), Survival and growth with a liability: Optimal portfolio strategies in continuous time, {\it Mathematics of Operations Research}, 22 (2): 468-493.

\smallskip \noindent \hangindent 20 pt Browne, S. (1999a), Beating a moving target: Optimal portfolio strategies for outperforming a stochastic benchmark, {\it Finance and Stochastics}, 3: 275-294.

\smallskip \noindent \hangindent 20 pt Browne, S. (1999b), The risk and rewards of minimizing shortfall probability, {\it Journal of Portfolio Management}, 25 (4): 76-85.

\smallskip \noindent \hangindent 20 pt Davis, M. H. A. and A. R. Norman (1990) Portfolio selection with transaction costs, {\it Mathematics of Operations Research}, 15: 676-713.

\smallskip \noindent \hangindent 20 pt Dhaene, J., S. Vanduffel, M. J. Goovaerts, R. Kaas, and D. Vyncke (2005), Comonotonic approximations for optimal portfolio selection problems, {\it Journal of Risk and Insurance}, 72 (2): 253-300.

\smallskip \noindent \hangindent 20 pt Duffie, D., W. Fleming, H. M. Soner, and T. Zariphopoulou (1997), Hedging in incomplete markets with HARA utility, {\it Journal of Economic Dynamics and Control}, 21: 753-782.

\smallskip \noindent \hangindent = 20 pt Duffie, D. and T. Zariphopoulou (1993), Optimal investment with undiversifiable income risk, {\it Mathematical Finance}, 3: 135-148.

\smallskip \noindent \hangindent = 20 pt Fleming, W. H. and T. Zariphopoulou (1991), An optimal investment/consumption model with borrowing, {\it Mathematics of Operations Research}, 16 (4): 802-822.

\smallskip \noindent \hangindent = 20 pt Hipp, C. and M. Plum (2000), Optimal investment for insurers, {\it Insurance: Mathematics and Economics}, 27: 215-228.

\smallskip \noindent \hangindent = 20 pt Hipp, C. and M. Taksar (2000), Stochastic control for optimal new business, {\it Insurance: Mathematics and Economics}, 26: 185-192.

\smallskip \noindent \hangindent = 20 pt Karatzas, I. and S. E. Shreve (1991), {\it Brownian Motion and Stochastic Calculus}, second edition, Springer-Verlag, New York.

\smallskip \noindent \hangindent = 20 pt Karatzas, I. and S. Shreve (1998), {\it Methods of Mathematical Finance}, Springer-Verlag, New York.

\smallskip \noindent \hangindent 20 pt Koo, H. K. (1998), Consumption and portfolio selection with labor income: A continuous time approach, {\it Mathematical Finance}, 8: 49-65.

\smallskip \noindent \hangindent = 20 pt Merton, R. C. (1992), {\it Continuous-Time Finance}, revised edition, Blackwell Publishers, Cambridge, Massachusetts.

\smallskip \noindent \hangindent = 20 pt Milevsky, M. A., K. Ho, and C. Robinson (1997), Asset allocation via the conditional first exit time or how to avoid outliving your money, {\it Review of Quantitative Finance and Accounting}, 9 (1): 53-70.

\smallskip \noindent \hangindent 20 pt Milevsky, M. A., K. S. Moore, and V. R. Young (2004), Optimal asset allocation and ruin-minimization annuitization strategies: The fixed consumption case, working paper, Department of Mathematics, University of Michigan.

\smallskip \noindent \hangindent = 20 pt Milevsky, M. A. and C. Robinson (2000), Self-annuitization and ruin in retirement, with discussion, {\it North American Actuarial Journal}, 4 (4): 112-129.

\smallskip \noindent \hangindent = 20 pt Olivieri, A. and E. Pitacco (2003), Solvency requirements for pension annuities, {\it Journal of Pension Economics and Finance}, 2 (2): 127-157.

\smallskip \noindent \hangindent 20 pt Parikh, A. N. (2003), The evolving U.S. retirement system, {\it The Actuary}, March: 2-6.

\smallskip \noindent \hangindent 20 pt Pestien, V. C. and W. D. Sudderth (1985), Continuous-time red and black: How to control a diffusion to a goal, {\it Mathematics of Operations Research}, 10 (4): 599-611.

\smallskip \noindent \hangindent 20 pt Rockafellar, R. T. (1970), {\it Convex Analysis}, Princeton University Press, Princeton, New Jersey.

\smallskip \noindent \hangindent = 20 pt Schmidli, H. (2001), Optimal proportional reinsurance policies in a dynamic setting, {\it Scandinavian Actuarial Journal}, 2001 (1): 55-68.

\smallskip \noindent \hangindent 20 pt Shreve, S. E. and H. M. Soner (1994), Optimal investment and consumption with transaction costs, {\it Annals of Applied Probability}, 4 (3): 206-236.

\smallskip \noindent \hangindent = 20 pt VanDerhei, J. and C. Copeland (2003), Can America afford tomorrow's retirees: Results from the EBRI-ERF Retirement Security Projection Model, Issue Brief of the Employee Benefit Research Institute, www.ebri.org, November 2003.

\smallskip \noindent \hangindent = 20 pt Vanduffel, S., J. Dhaene, and M. Goovaerts (2005), On the evaluation of `saving-consump-tion' plans, {\it Journal of Pension Economics and Finance}, 4 (1): 17-30.

\smallskip \noindent \hangindent = 20 pt Young, V. R. (2004), Optimal investment strategy to minimize the probability of lifetime ruin, to appear, {\it North American Actuarial Journal}.

\smallskip \noindent \hangindent 20 pt Zariphopoulou, T. (1992), Investment/consumption models with transaction costs and Markov-chain parameters, {\it SIAM Journal on Control and Optimization}, 30: 613-636.

\smallskip \noindent \hangindent = 20 pt Zariphopoulou, T. (1994), Consumption-investment models with constraints, {\it SIAM Journal on Control and Optimization}, 32 (1): 59-85.

\smallskip \noindent \hangindent 20 pt Zariphopoulou, T. (1999), Transaction costs in portfolio management and derivative pricing, {\it Introduction to Mathematical Finance}, D. C. Heath and G. Swindle (editors), American Mathematical Society, Providence, RI. {\it Proceedings of Symposia in Applied Mathematics}, 57: 101-163.

\smallskip \noindent \hangindent 20 pt Zariphopoulou, T. (2001), Stochastic control methods in asset pricing, {\it Handbook of \break Stochastic Analysis and Applications}, D. Kannan and V. Lakshmikantham (editors), Marcel Dekker, New York.

\vfill
\eject

\centerline{\vbox{\hsize=4.5in\psfig{figure=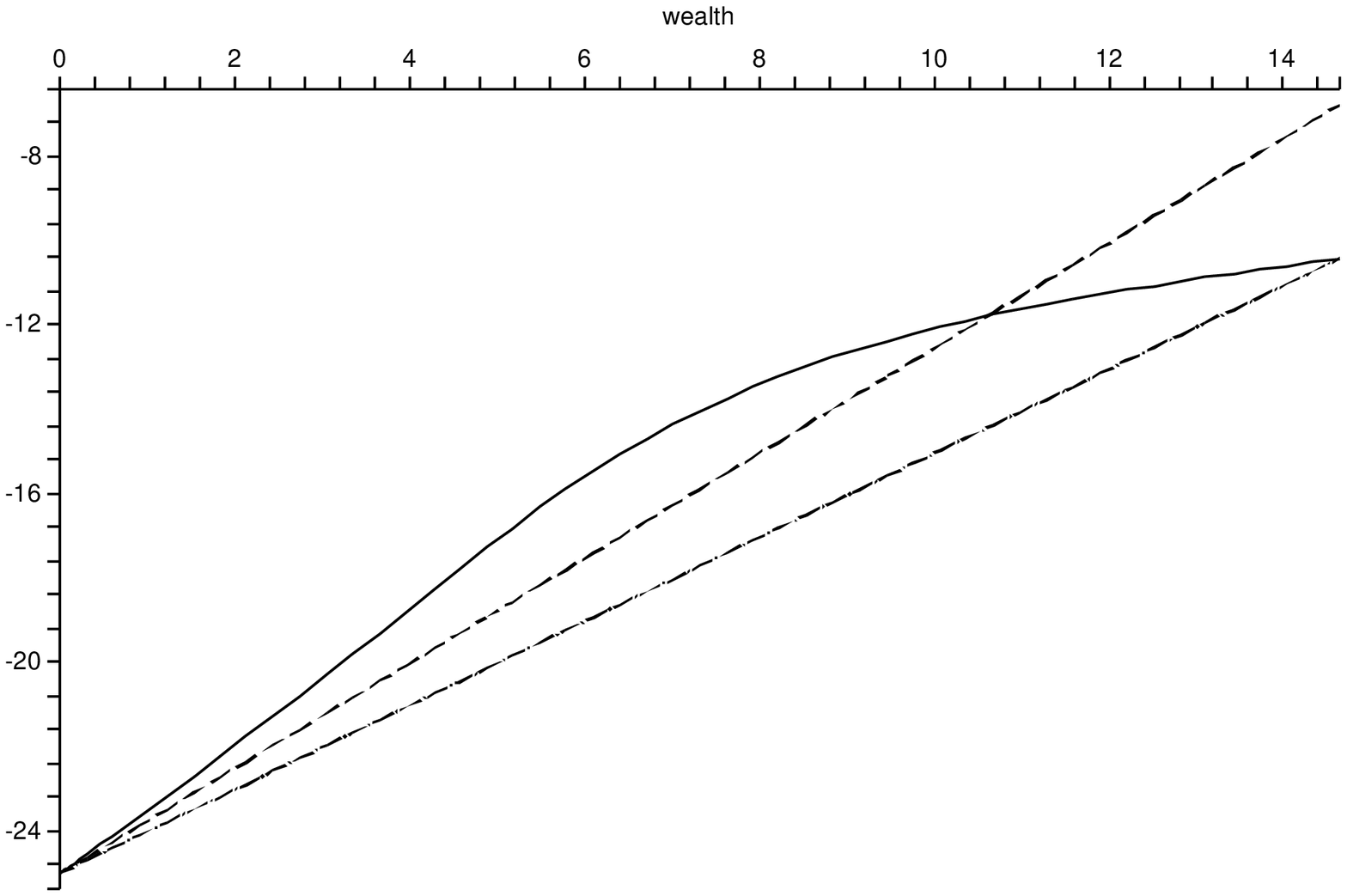,width=4.5in}
   \hfill \break   
   Figure 1:  Graph of $y$, $z$, and $z_b$.  The solid line corresponds to $y$; the dashed line to $z$; and the dotted line to $z_b$.}}

\centerline{\vbox{\hsize=4.5in\psfig{figure=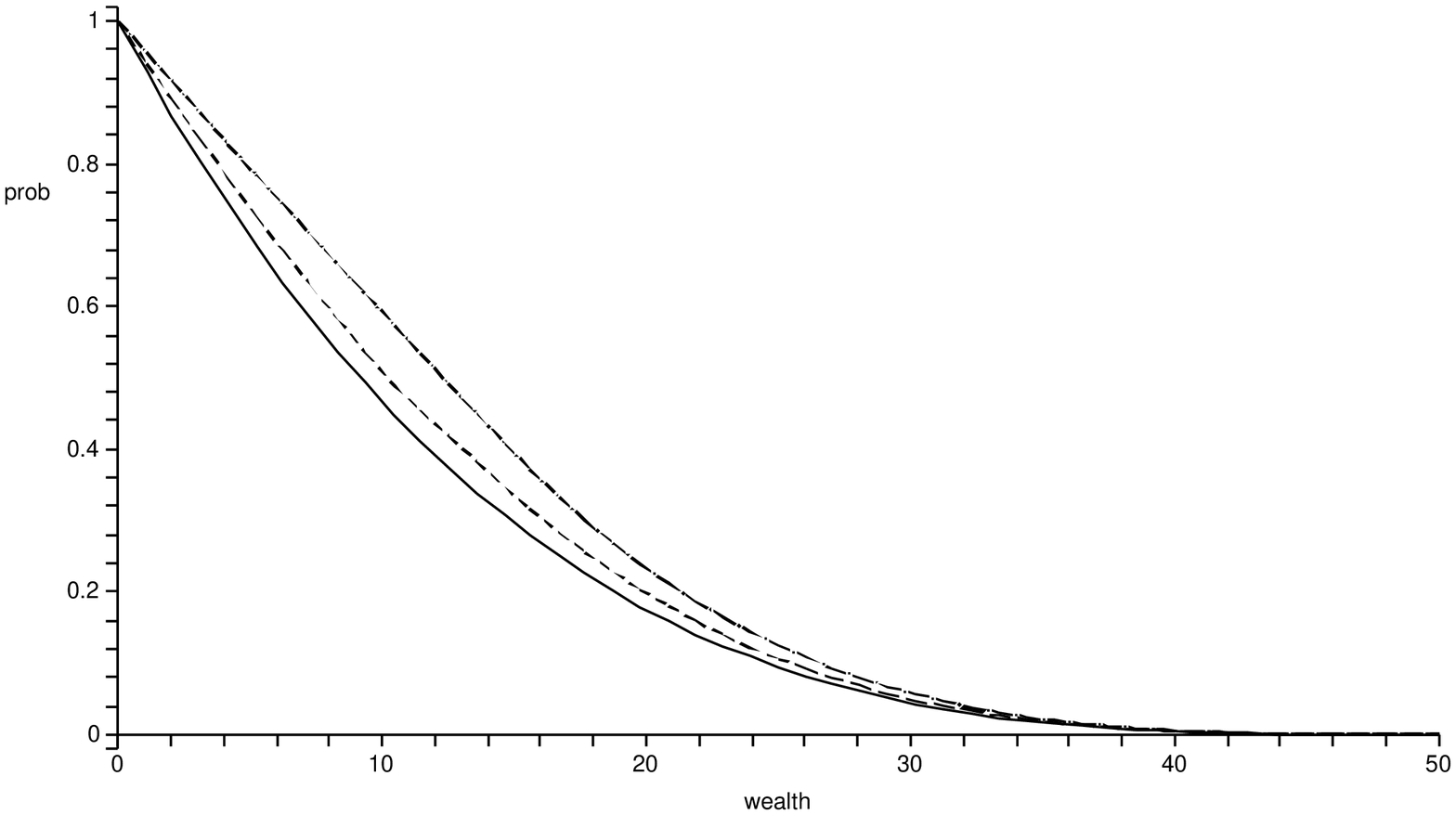,width=4.5in}
   \hfill \break   
   Figure 2:  Graph of $\psi$, $\psi_0$, and $\psi_b$.  The solid line corresponds to $\psi$; the dashed line to $\psi_0$; and the dotted line to $\psi_b$.}}

\vfill
\eject

\centerline{\vbox{\hsize=4.5in\psfig{figure=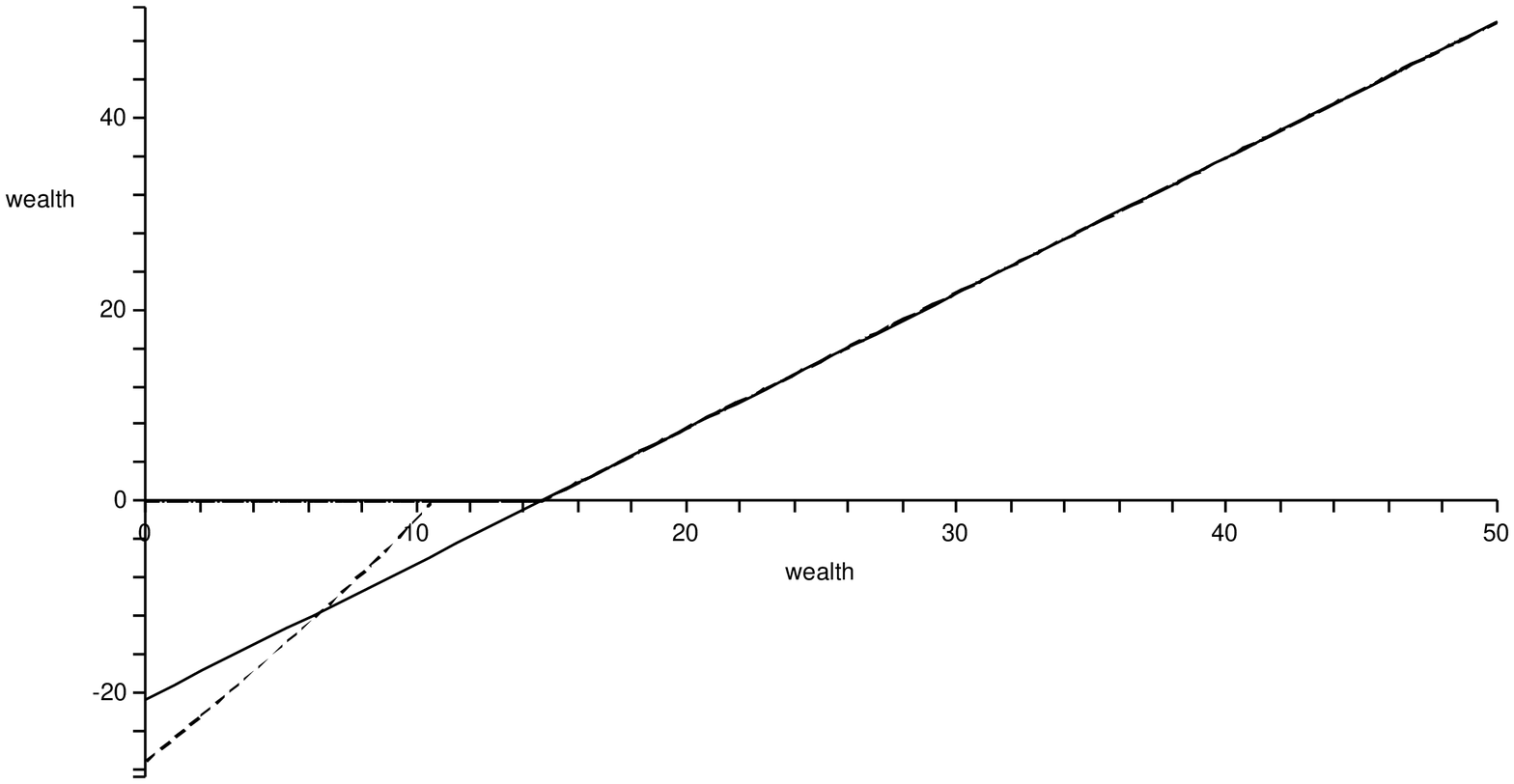,width=4.5in} 
\hfill \break   
Figure 3:  Graph of $\zeta$, $\zeta_0$, and $\zeta_b$.  The solid line corresponds to $\zeta$; the dashed line to $\zeta_0$; and the dotted line to $\zeta_b$.}}

\centerline{\vbox{\hsize=4.5in\psfig{figure=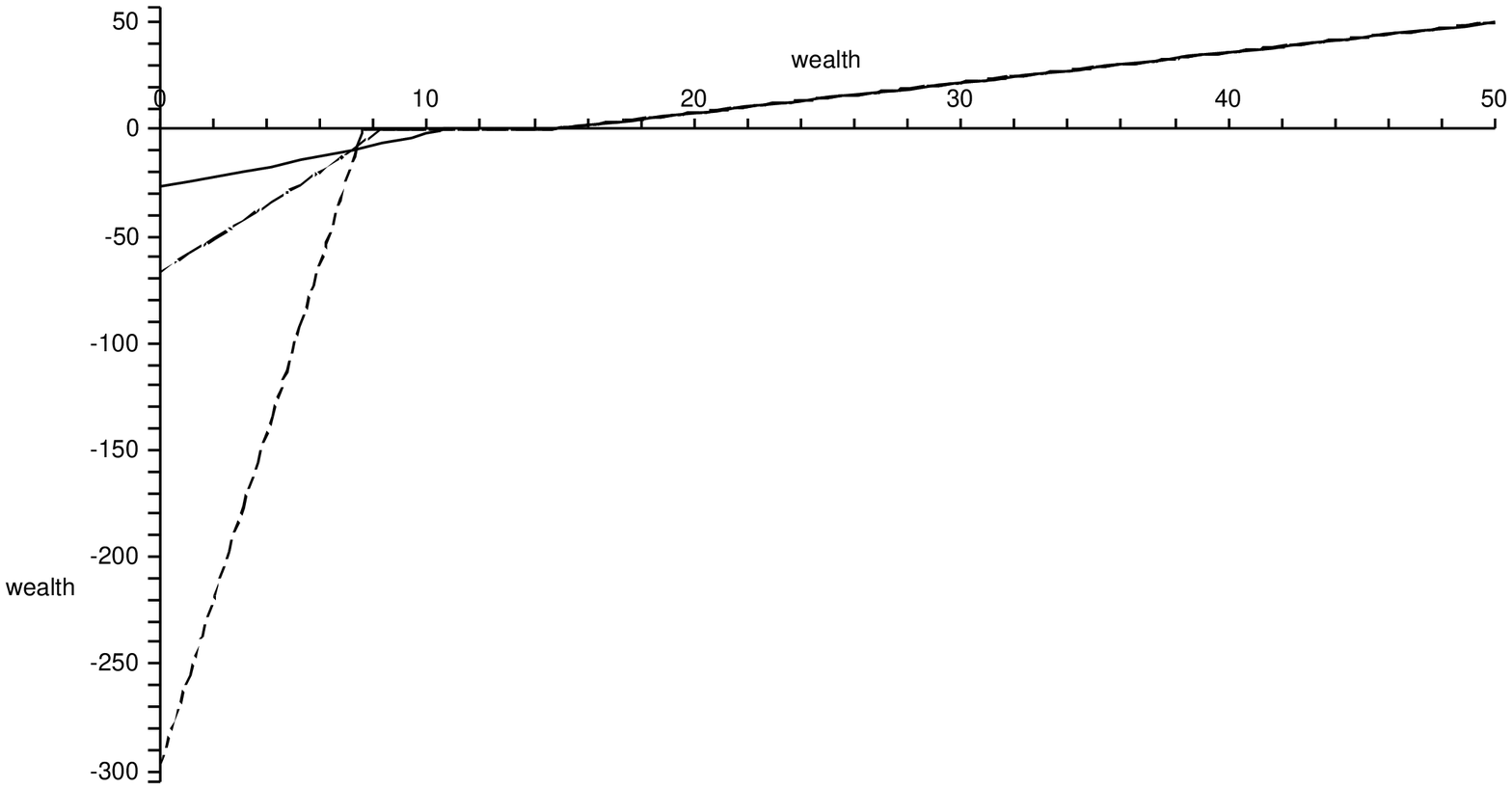,width=4.5in} 
\hfill \break 
Figure 4:  Graph of $\zeta_b$ for various values of $b$.  The solid line corresponds to $b = 0.04$; the dashed line to $b = 0.055$; and the dotted line to $b = 0.59$.}}

\bye